\documentclass[oneside,english]{amsart}
\usepackage[T1]{fontenc}
\usepackage[latin1]{inputenc}
\usepackage{prettyref}
\usepackage{graphicx}
\usepackage{amssymb}

\makeatletter

\providecommand{\tabularnewline}{\\}

 \theoremstyle{plain}    
 \newtheorem{thm}{Theorem}[section]
 \numberwithin{equation}{section} 
 \numberwithin{figure}{section} 
 \theoremstyle{plain}
 \theoremstyle{plain}    
 \newtheorem{cor}[thm]{Corollary} 
 \theoremstyle{plain}    
 \newtheorem{lem}[thm]{Lemma} 
 \theoremstyle{plain}    
 \newtheorem{prop}[thm]{Proposition} 
 \theoremstyle{definition}
 \newtheorem{defn}[thm]{Definition}
 \theoremstyle{definition}
  \newtheorem*{example*}{Example}
 \theoremstyle{remark}    
 \newtheorem{claim}[thm]{Claim}
 \theoremstyle{remark}    
 \newtheorem*{claim*}{Claim}

\usepackage{fullpage}

\newrefformat{pro}{Proposition #1}
\newrefformat{cor}{Corollary #1}
\newrefformat{def}{Definition #1}
\newrefformat{cla}{Claim #1}

\newcommand{\POlong}{\textsc{orderable}}
\newcommand{\PO}{\textsc{po}}

\newcommand{\PRlong}{\textsc{reachable}}
\newcommand{\PR}{\textsc{pr}}
\newcommand{\NPR}{\textsc{npr}}
\newcommand{\PClong}{\textsc{coverable}}
\newcommand{\PC}{\textsc{pc}}
\newcommand{\PAlong}{\textsc{annihilation}}
\newcommand{\PA}{\PAlong}
\newcommand{\TSAT}{\textsc{3sat}}
\newcommand{\RTSAT}{\textsc{r3sat}}

\newcommand{\OPN}{\textsc{opn}}
\newcommand{\OPNlong}{\textsc{optimal-pebbling-number}}
\newcommand{\DPR}{\textsc{dpr}}
\newcommand{\PN}{\textsc{pn}}
\newcommand{\PNlong}{\textsc{pebbling-number}}
\newcommand{\RPNlong}{\textsc{r-pebbling-number}}
\newcommand{\RPN}{\textsc{rpn}}
\newcommand{\AETSAT}{\forall\exists\TSAT}
\newcommand{\RAETSAT}{\textsc{r}\AETSAT}

\newcommand{\ptc}{\Pi_2^{\mathrm{P}}}
\newcommand{\balance}{\mathrm{balance}}
\newcommand{\starg}{\mathrm{star}}

\newcommand{\comp}{\mathrm{comp}}
\newcommand{\piopt}{\widehat{\pi}}

\newcommand{\poly}{\mathrm{poly}}

\usepackage{babel}
\makeatother
\begin{document}

\title{The Complexity of Graph Pebbling}

\author{Kevin Milans}

\address{Kevin Milans\\
Computer Science Department\\
University of Illinois, Urbana, IL, 61801}

\email{milans@uiuc.edu}

\author{Bryan Clark}

\address{Bryan Clark\\
Physics Department\\
University of Illinois, Urbana, IL, 61801}

\email{bkclark@uiuc.edu}

\date{\today}

\begin{abstract}
In a graph $G$ whose vertices contain pebbles, a pebbling move $uv$
removes two pebbles from $u$ and adds one pebble to a neighbor $v$
of $u$. The optimal pebbling number $\piopt(G)$ is the minimum $k$
such that \emph{there exists} a distribution of $k$ pebbles to $G$
so that for any target vertex $r$ in $G$, there is a sequence of
pebbling moves which places a pebble on $r$. The pebbling number
$\pi(G)$ is the minimum $k$ such that \emph{for all} distributions
of $k$ pebbles to $G$ and for any target vertex $r$, there is a
sequence of pebbling moves which places a pebble on $r$. 

We explore the computational complexity of computing $\piopt(G)$
and $\pi(G)$. In particular, we show that deciding whether $\piopt(G)\leq k$
is NP-complete and deciding whether $\pi(G)\leq k$ is $\ptc$-complete.
Additionally, we provide a characterization of when an unordered set
of pebbling moves can be ordered to form a valid sequence of pebbling
moves.
\end{abstract}
\maketitle

\section{Introduction}

Let $G$ be a simple, undirected graph and let $p\,:\, V(G)\to\mathbb{N}\cup\{0\}$
be a distribution of pebbles to the vertices of $G$. We refer to
the total number of pebbles $\sum_{v}p(v)$ as the \emph{size} of
$p$, denoted by $\left|p\right|$. A pebbling move $uv$ \emph{}consists
of removing two pebbles from a vertex $u$ with $p(u)\geq2$ and placing
one pebble on a neighbor $v$ of $u$. After completing a pebbling
move $uv$, we are left with a new distribution of pebbles, which
we denote by $p_{uv}$. Similarly, if $\sigma=u_{1}v_{1},\ldots,u_{k}v_{k}$
is a sequence of pebbling moves, denote by $p_{\sigma}$ the distribution
of pebbles that results from making the pebbling moves specified by
$\sigma$. Although graph pebbling was originally developed to simplify
a result in number theory (F.R.K. Chung provides the history \cite{history}),
it has since become an object of study in its own right. G.H. Hurlbert
presents a detailed survey of early graph pebbling results \cite{survey-paper}. 

\begin{description}
\item [Notational~Conventions]We use $G$ and $H$ to refer to simple,
undirected graphs. We use $D$ and $E$ to refer to directed graphs
with multiple edges and loops. If $v$ is a vertex in a directed multigraph,
we denote the indegree (resp. outdegree) of $v$ by $d^{+}(v)$ (resp.
$d^{-}(v)$). We write $n(G)$ (resp. $e(G)$) for the number of vertices
(resp. edges) in $G$. Similarly, we use $V(G)$ (resp. $E(G)$) to
refer to the vertex set (resp. edge set) of $G$. We write $d_{G}(u,v)$
(or $d(u,v)$ when $G$ is clear from context) for the length of the
shortest $uv$-path in $G$.\\
If $p$ and $q$ are pebble distributions on a graph $G$, we say
that $p\geq q$ if $p(v)\geq q(v)$ for each vertex $v$ in $G$.
\end{description}
Given a graph $G$ with a pebble distribution $p$, we say that a
vertex $r$ in $G$ is \emph{reachable} if there is a sequence of
pebbling moves which places a pebble on $r$. Note that whenever $p(r)>0$,
$r$ is trivially reachable. The notion of reachability is fundamental
to graph pebbling; most of our decision problems involve questions
of reachability. We call the problem of deciding (given $G$, $p$,
and $r$) whether $r$ is reachable $\PRlong$. In section 3, we establish
that $\PRlong$ is NP-complete, a result obtained simultaneously and
independently by N.G. Watson \cite{watson}.

Given a graph $G$ and a target vertex $r$, the $r$-\emph{pebbling
number} of $G$, denoted $\pi(G,r)$, is the minimum $k$ such that
$r$ is reachable under every pebble distribution of size $k$. Similarly,
the \emph{pebbling number} of $G$, denoted $\pi(G)$, is the minimum
$k$ such that every vertex in $G$ is reachable under every pebble
distribution of size $k$. For a connected graph $G$, a pigeonhole
argument quickly establishes that such a $k$ exists and so $\pi(G)$
is well defined (see \prettyref{pro:PN-trivial-bound}). We call the
problem of deciding whether $\pi(G,r)\leq k$ (resp. $\pi(G)\leq k$)
$\RPNlong$ (resp. $\PNlong$). In section 5, we establish that both
decision problems are $\ptc$-complete, meaning that these problems
are complete for the class of problems computable in polynomial time
by a coNP machine equipped with an oracle for an NP-complete language.
Consequently, these decision problems are both NP-hard and coNP-hard.
It follows that $\RPNlong$ and $\PNlong$ are neither in NP nor in
coNP unless $\textrm{NP}=\textrm{coNP}$. N.G. Watson simultaneously
and independently established that $\RPNlong$ is coNP-hard \cite{watson}.

Observe that if we fix some $r$ in $G$ and put one pebble on every
other vertex, $r$ is not reachable. It follows that $\pi(G)\geq n(G)$.
It is natural to wonder which graphs achieve equality in $\pi(G)=n(G)$.
Although no characterization of such graphs is known, a growing body
of results provide conditions that are necessary or sufficient to
imply $\pi(G)=n(G)$. Recall that $G$ is $k$-connected if $n(G)\geq k+1$
and for every set $S$ of at most $k-1$ vertices, $G-S$ is connected.
In \cite{diameter-two-graphs}, it is shown that if $G$ is $3$-connected
and has diameter $2$, then $\pi(G)=n(G)$. Consequently, the probability
that a random graph on $n$ vertices satisfies $\pi(G)=n(G)$ approaches
$1$ as $n$ grows. Furthermore, in \cite{large-connectivity} it
is shown that if $G$ has diameter $d$ and is $(2^{2d+3})$-connected,
then $\pi(G)=n(G)$. On the other hand, if $G$ contains a cut vertex,
then $\pi(G)>n(G)$. Indeed, suppose $v$ is a cut vertex in $G$
and let $u$ and $w$ be vertices in separate components of $G-v$.
If we put three pebbles on $u$, zero pebbles on $v$ and $w$, and
one pebble on every other vertex, then it is not possible to place
a pebble on $w$. 

The \emph{optimal pebbling number} of $G$, denoted $\piopt(G)$,
is the minimum $k$ such that each vertex is reachable under some
distribution of size $k$. We call the problem of deciding whether
$\piopt(G)\leq k$ $\OPNlong$. In section 4, we establish that $\OPNlong$
is NP-complete. 

It is immediate that $\piopt(G)\leq n(G)$. However, this bound is
not tight for connected graphs. As shown in \cite{shameless-plug},
if $G$ is connected, then $\piopt(G)\leq\left\lceil 2n(G)/3\right\rceil $.
Equality is achieved by the path \cite{shameless-plug,PSV} and the
cycle \cite{shameless-plug}. It is an open problem to characterize
which graphs achieve equality.

Given $G$ and distributions $p$ and $q$, we say that $p$ \emph{covers}
$q$ if there exists a sequence of pebbling moves $\sigma$ such that
$p_{\sigma}\geq q$. The \emph{unit distribution} assigns one pebble
to each vertex in $G$. We call the problem of deciding whether $p$
covers the unit distribution $\PClong$. In section 3, we establish
that $\PClong$ is NP-complete; this result was obtained simultaneously
and independently by N.G. Watson \cite{watson}. 

Although most of the problems we study are computationally difficult,
there are some interesting pebbling problems that are tractable. A
pebble distribution $q$ is \emph{positive} if $q$ assigns at least
one pebble to every vertex. A distribution $p$ is \emph{simple} if
it assigns zero pebbles to all but one vertex. The \emph{$q$-cover
pebbling number} of $G$, denoted $\gamma_{q}(G)$, is the minimum
$k$ such that every distribution of size $k$ covers $q$. The cover
pebbling theorem states that for any positive distribution $q$, there
is a simple distribution $p$ of size $\gamma_{q}(G)-1$ such that
$p$ does not cover $q$ \cite{cover-thm-1,cover-thm-2}. As a consequence,
given $G$ and a positive distribution $q$, one easily computes $\gamma_{q}(G)$
in polynomial time. In the special case that $q$ is the unit distribution,
we simply write $\gamma(G)$ for $\gamma_{q}(G)$.

\begin{figure}
\begin{center}\begin{tabular}{|c|c|c|c|}
\hline 
Short Name&
Full Name&
Description&
Complexity\tabularnewline
\hline
\hline 
$\PN$&
$\PNlong$&
Given $G$,$k$: is $\pi(G)\leq k$?&
$\ptc$-complete\tabularnewline
\hline 
$\RPN$&
$\RPNlong$&
Given $G,k,r$: is $\pi(G,r)\leq k$?&
$\ptc$-complete\tabularnewline
\hline 
$\OPN$&
$\OPNlong$&
Given $G,k$: is $\piopt(G)\leq k$?&
NP-complete\tabularnewline
\hline 
$\PR$&
$\PRlong$&
Given $G,p,r$: is $r$ reachable?&
NP-complete\tabularnewline
\hline 
$\PC$&
$\PClong$&
Given $G,p$: does $p$ cover the unit distribution?&
NP-complete\tabularnewline
\hline
\end{tabular}\end{center}

\caption{A summary of the decision problems considered in this paper\label{fig:results-summary}}
\end{figure}

\begin{description}
\item [Overview]In section 2, we develop a characterization of when unordered
sets of pebbling moves may ordered in a way that yields a valid sequence
of pebbling moves. In section 3, we present results on the complexity
of $\PRlong$ and $\PClong$. We also observe that a simple greedy
strategy solves $\PRlong$ whenever $G$ is a tree. Section 3 uses
some results from section 2. In section 4, we present our results
on the complexity of computing the optimal pebbling number. Section
4 uses some results from sections 2 and 3. In section 5, we present
our results on the complexity of computing the ($r$-)pebbling number.
This section uses some results from sections 2 and 3; it is generally
independent of section 4. In section 6, we present our conclusions. 
\end{description}
Let us consider a simple example. Suppose we are given a graph $H$
with a distribution of pebbles, and we wish to determine if there
is a sequence of pebbling moves which ends with only one pebble left
in the entire graph. We call this problem $\PAlong$. It is not difficult
to see that $\PAlong$ is NP-hard. Indeed, a reduction from $\textsc{hamiltonian-path}$
is almost immediate. Specifically, to decide if $G$ has a Hamiltonian
path, we may construct $H$ from $G$ by introducing a new vertex
$v$ which is adjacent to each vertex in $G$. We place two pebbles
on $v$ and one pebble on every other vertex in $H$. It is clear
that $G$ has a Hamiltonian path if and only if there is a sequence
of pebbling moves which results in only one pebble in $H$. 

What is less clear is that $\PA$ is in NP. If $\sigma$ is a sequence
of pebbling moves in $G$ under $p$ which results in only one pebble
left in $G$, then the length of $\sigma$ is $\left|p\right|-1$,
which may be exponentially large in the number of bits needed to represent
$G$ and $p$. Hence, $\sigma$ may be too large to serve as a certificate
for membership in $\PA$. However, as we will see, the order of the
moves in $\sigma$ is insignificant. In fact, if we are merely told
how many times $\sigma$ pebbles along each direction in every edge
in $G$, then we can quickly verify the existence of $\sigma$.

\section{Pebble Orderability}

Many questions in graph pebbling concern the existence of a sequence
of pebbling moves with certain properties. There is a natural temptation
to search for such sequences directly, by deciding which pebbling
move to make first, which to make second, and so forth. In this section,
we develop tools that allow us more flexibility in constructing sequences
of pebbling moves. In particular, our goal is to worry only about
which moves we should make, and not the order in which to make them.

We define the \emph{signature} of a sequence of pebbling moves $\sigma$
in a graph $G$ to be the directed multigraph on vertex set $V(G)$
where the multiplicity of an edge $uv$ is the number of times $\sigma$
pebbles from $u$ to $v$. We say that a digraph $D$ is \emph{orderable}
under a pebble distribution $p$ if some ordering of $E(D)$ is a
valid sequence of pebbling moves. We characterize when $D$ is orderable
under $p$. We call the problem of testing whether $D$ is orderable
under $p$ $\POlong$, or $\PO$ for short.

As it turns out, two conditions which are necessary for $D$ to be
orderable are also sufficient. Suppose that $D$ is orderable and
consider a vertex $v$. We note that $v$ begins with $p(v)$ pebbles,
$D$ pledges that $v$ will receive $d_{D}^{-}(v)$ pebbles from pebbling
moves into $v$, and $D$ requests $d_{D}^{+}(v)$ pebbling moves
out of $v$. Because each pebbling move out of $v$ costs two pebbles,
it is clear that $p(v)+d_{D}^{-}(v)$ is at least $2d_{D}^{+}(v)$.
This leads us to define the \emph{balance} of a vertex $v$ as \[
\balance(D,p,v)=p(v)+d_{D}^{-}(v)-2d_{D}^{+}(v).\]
 The balance of $v$ is simply the number of pebbles that remain on
$v$ after executing any sequence of pebbling moves whose signature
is $D$; that is, for any $\sigma$ whose signature is $D$, we have
that $p_{\sigma}(v)=\balance(D,p,v)$. 

If $D$ is orderable under $p$, then the balance of each vertex must
be nonnegative. We call this condition the \emph{balance condition}.
The balance condition alone is not sufficient: if $D$ is a directed
cycle and each vertex has one pebble, then the balance of each vertex
is zero but we cannot make any pebbling moves, and so $D$ is not
orderable. However, as was implicitly observed in \cite{Moews}, if
$D$ is acyclic, then the balance condition is sufficient. 

\begin{thm}
[Acyclic Orderability Characterization]\cite{Moews} \label{thm:acyclic-feasibility}
If $D$ is an acyclic digraph with distribution $p$, then $D$ is
orderable if and only if the balance condition is satisfied.
\end{thm}
\begin{proof}
We have observed that the balance condition is necessary. Conversely,
if the balance condition is satisfied, then we obtain a sequence of
pebbling moves $\sigma$ whose signature is $D$ by iteratively selecting
a source $u$ in $D$, making all pebbling moves out of $u$, and
deleting $u$ from $D$. 
\end{proof}
Despite the simplicity of the acyclic orderability characterization,
we are already able to obtain one of our most useful corollaries.
It makes precise our intuition that if we are trying to place pebbles
on a target vertex $r$, it is never advantageous to pebble around
in a cycle. Our proof is somewhat shorter than previous proofs.

\begin{cor}
[No Cycle Lemma]\cite{cover-pebbling,Moews} \label{cor:cycles-can-be-removed}Suppose
$D$ is orderable under $p$. There exists an acyclic $D'\subseteq D$
such that $D'$ is orderable and $\balance(D',p,v)\geq\balance(D,p,v)$
for all $v$. 
\end{cor}
\begin{proof}
Let $D'$ be a digraph obtained by iteratively removing cycles from
$D$ until no cycles remain. Observe that removing a cycle $C$ does
not change the balance of vertices outside of $C$ and increases the
balance of vertices in $C$ by one. It follows that $\balance(D',p,v)\geq\balance(D,p,v)\geq0$
for all $v$. Hence, $D'$ is acyclic and satisfies the balance condition.
By the acyclic orderability characterization, $D'$ is orderable. 
\end{proof}
In most contexts, if a sequence $\sigma$ of pebbling moves satisfies
certain criterion, then so will any sequence $\sigma'$ provided that
$p_{\sigma'}\geq p_{\sigma}$. As we have seen, in these situations,
we are able to restrict our attention to sequences of pebbling moves
whose signatures are acyclic. Indeed, all of our major results fall
into this category and therefore only require the orderability characterization
for acyclic digraphs. 

Nevertheless, one may wish to study the existence of sequences of
pebbling moves which purposefully remove pebbles from the graph, as
in the $\PAlong$ decision problem. Let us return to our orderability
characterization for arbitrary $D$. As we have seen, in general the
balance condition is not sufficient. However, as we show in our next
lemma, a directed cycle with one pebble on each vertex is the only
minimal, nontrivial situation which satisfies the balance condition
and does not allow us to make any pebbling moves. 

\begin{lem}
\label{lem:unique-stop-configuration} Suppose that $D$ with distribution
$p$ satisfies the balance condition, $D$ is connected, and $e(D)\geq1$.
If we cannot make any pebbling move described by an edge in $D$,
then $D$ is a directed cycle and each vertex has exactly one pebble. 
\end{lem}
\begin{proof}
Observe that $D$ does not have any source vertices. Indeed, if $v$
were a source, then the balance condition implies that $v$ has enough
pebbles to make all pebbling moves out of $v$ requested by $D$.
Therefore $v$ must have outdegree zero, and so $v$ is an isolated,
loopless vertex, which contradicts that $D$ is connected and contains
an edge. 

Let $e=E(D)$, let $X\subseteq V(D)$ be the set of all sinks, let
$Y=V(D)-X$ be the set of all nonsinks, let $k$ be the number of
edges with sources in $Y$ and sinks in $X$, and let $z$ be the
number of nonsinks that have exactly one pebble. Note that $e=\sum_{v}d^{-}(v)=k+\sum_{v\in Y}d^{-}(v)$
and $e=\sum_{v}d^{+}(v)=\sum_{v\in Y}d^{+}(v)$. Furthermore, for
each $v\in Y$, we have that $p(v)\leq1$; otherwise, $p(v)\geq2$
and $v$ has outdegree at least one, contradicting that there are
no pebbling moves available. It follows that $z=\sum_{v\in Y}p\left(v\right)$.
Adding the inequality $\balance(D,p,v)\geq0$ over all $v\in Y$,
we obtain \[
\sum_{v\in Y}d^{-}(v)+\sum_{v\in Y}p\left(v\right)\geq2\sum_{v\in Y}d^{+}(v),\]
or equivalently $e-k+z\geq2e$, and so $e+k\leq z$. Because $D$
has no sources, every vertex has indegree at least one and so $e\geq n$.
Therefore $n\leq e\leq e+k\leq z\leq n$. It follows that $n=e=z$,
so that every vertex in $D$ is neither a sink nor a source and has
exactly one pebble. Furthermore, because $e=n$, each vertex in $D$
has indegree and outdegree exactly one. It follows that $D$ is a
directed cycle.
\end{proof}
Of course, any sequence of pebbling moves leaves a pebble somewhere
in the graph; therefore if $D$ contains an edge and $D$ is orderable,
then $\balance(D,p,v)\geq1$ for some vertex $v$. In fact, a slight
generalization of this observation will serve as our second necessary
condition. To develop this condition, we first recall the component
digraph.

Let $D$ be a directed multigraph. A strongly connected component
$A$ is \emph{trivial} if $A$ consists of a single vertex with indegree
and outdegree zero. Define $\comp(D)$, the \emph{component digraph}
of $D$, to be the digraph obtained by contracting each strongly connected
component of $D$ to a single vertex. 

Suppose that $D$ is orderable, and consider a sink $A$ in $\comp(D)$.
Because $A$ is a sink component, any pebbling move whose source is
in $A$ also has its sink in $A$; it follows that unless $A$ is
trivial, then there must be some vertex $v$ in $A$ with $\balance(D,p,v)\geq1$.
We call the condition that every nontrivial sink in $\comp(D)$ contains
a vertex of positive balance the \emph{sink condition}. Note that
in the directed cycle example, each vertex has balance zero, and so
it fails the sink condition. 

As we now show, the balance condition together with the sink condition
are sufficient for $D$ to be orderable. We require a simple proposition. 

\begin{prop}
\label{pro:disconnect-D} If $D$ is a strongly connected digraph
and $D-uv$ is not strongly connected, then $\comp(D-uv)$ contains
a single sink $A$, $u$ is in $A$, and $v$ is not in $A$.
\end{prop}
\begin{thm}
[Orderability Characterization]\label{thm:feasibility-characterization}
$D$ is orderable under $p$ if and only if 
\begin{enumerate}
\item (balance condition) every vertex has nonnegative balance, and
\item (sink condition) every nontrivial sink $A$ in $\comp(D)$ contains
some vertex with balance at least one
\end{enumerate}
\end{thm}
\begin{proof}
We have observed that both conditions are necessary. We show that
$D$ is orderable under $p$ by induction on $e(D)$. If $e(D)=0$,
the statement is trivial. In the remaining cases, we assume that $D$
has at least one edge.

We consider the case that there is a source $v$ in $D$ with outdegree
at least one. Because $\balance(D,p,v)\geq0$, $v$ has enough pebbles
to make all the pebbling moves that $D$ requests out of $v$. Let
$\sigma$ be an arbitrary ordering of these moves and obtain $D'$
from $D$ by removing all edges whose source is $v$. We argue that
$D'$ is orderable under $p_{\sigma}$. It is clear that $D'$ under
$p_{\sigma}$ satisfies the balance condition. Observe that every
sink in $\comp(D')$ either consists of $v$ (and is therefore trivial)
or is a sink in $\comp(D)$. It follows that every nontrivial sink
in $\comp(D')$ is a nontrivial sink in $\comp(D)$ and hence contains
some vertex with balance at least one. By induction, $D'$ is orderable
under $p_{\sigma}$. In the remaining cases, we assume that every
source in $D$ is an isolated vertex.

Next, we consider the case where $\comp(D)$ contains a source $A$
with outdegree at least one. Let $uv$ be an edge from a vertex $u$
in $A$ to a vertex $v$ outside of $A$. We check that $A$ under
$p$ satisfies both the balance condition and the sink condition.
The balance condition follows from observing that $A$ is a source
in $\comp(D)$. Because $A$ is strongly connected and $\balance(A,p,u)\geq2$,
we have that $A$ satisfies the sink condition. By induction, there
is an ordering $\sigma$ of $E(A)$ which is a valid sequence of pebbling
moves. We argue that $D-E(A)$ is orderable under $p_{\sigma}$. It
is clear that $D-E(A)$ under $p_{\sigma}$ satisfies the balance
condition. Because every nontrivial sink in $\comp(D-E(A))$ is a
nontrivial sink in $\comp(D)$, $D-E(A)$ satisfies the sink condition.
Because every source in $D$ is an isolated vertex, it must be that
there is some edge $e$ in $D$ whose sink is $u$; this edge $e$
is contained in $A$. Therefore $D-E(A)$ contains fewer edges in
$D$, so that the inductive hypothesis implies that $D-E(A)$ is orderable
under $p_{\sigma}$. In the remaining cases we assume that every source
in $\comp(D)$ is an isolated vertex in $\comp(D)$.

Because $\comp(D)$ is acyclic and every source in $\comp(D)$ is
an isolated vertex in $\comp(D)$, it follows that $D$ consists of
disjoint, strongly connected components. Because $D$ is orderable
if and only if each component of $D$ is orderable, we assume without
loss of generality that $D$ is a single, strongly connected component.
If we can make a pebbling move $uv$ which leaves $D-uv$ strongly
connected, then it is clear that $D-uv$ under $p_{uv}$ satisfies
both conditions and so $D$ is orderable.

It remains to consider the case that every possible pebbling move
results in a digraph which is no longer strongly connected. By \prettyref{lem:unique-stop-configuration},
we have that some pebbling move $uv$ is possible. 

First, suppose that $uv$ is the only edge out of $u$. Note that
because $D$ is strongly connected, $u$ must have indegree at least
one. Furthermore, because $uv$ is a valid pebbling move, we have
$p\left(u\right)\geq2$. It follows that $\balance(D,p,u)\geq1$.
It is clear that $D-uv$ under $p_{uv}$ satisfies the balance condition;
together with \prettyref{pro:disconnect-D}, we have that it also
satisfies the sink condition. By induction $D-uv$ is orderable under
$p_{uv}$.

Otherwise, let $uw\in E\left(D\right)$, $w\neq v$. Let $z$ be a
vertex in $D$ with $\balance(D,p,z)\geq1$ (we allow $z\in\left\{ u,w,v\right\} $),
let $P$ be a $uz$-path, and let $Q$ be a $zu$-path. Observe that
$uv\not\in P$ or $uw\not\in P$. In the former case, $u$ and $z$
are in the same strongly connected component in $D-uv$; in the latter
case, $u$ and $z$ are in the same strongly connected component in
$D-uw$. Together with \prettyref{pro:disconnect-D} we have that
either $D-uv$ under $p_{uv}$ or $D-uw$ under $p_{uw}$ satisfies
both conditions. It follows that $D$ is orderable starting from $p$. 
\end{proof}
Observe that if $D$ is acyclic, then the sink condition is trivially
satisfied, and we recover the acyclic orderability characterization.
Our general orderability characterization yields a quick method for
checking whether $D$ is orderable, and so $\POlong$ is in P. As
a consequence, we see that $\PAlong$ is in NP. 

Before we conclude this section, we use our tools to prove some technical
lemmas which will be useful in later sections. 

\begin{lem}
\label{lem:large-outdegree-many-pebbles} Suppose $D$ is acyclic
and orderable under $p$. Then for any vertex $w$, there exists $D'\subseteq D$
such that $D'$ is orderable and \begin{eqnarray*}
\balance(D',p,w) & \geq & \balance(D,p,w)+2d_{D}^{+}(w)\\
 & \geq & p(w)+d_{D}^{-}(w).\end{eqnarray*}
 Additionally, if $d_{D}^{+}(w)>0$, then we may take $D'$ to be
a proper subgraph of $D$.
\end{lem}
\begin{proof}
Observe that if $uv$ is an edge in $D$ with $v$ a sink, then $D-uv$
satisfies the balance condition. Let $D'$ be a digraph obtained from
$D$ by iteratively deleting edges into sinks other than $w$ until
no such edges remain. Because $D'$ is acyclic and satisfies the balance
condition, the acyclic orderability characterization implies that
$D'$ is orderable. Observe that $w$ is a sink, or else $D'$ would
contain an edge $uv$ with $v\neq w$ a sink. Furthermore, every edge
into $w$ in $D$ remains in $D'$. It follows that $\balance(D',p,w)\geq\balance(D,p,w)+2d_{D}^{+}(w)$. 
\end{proof}
Often, we wish to explore the consequences of the existence of a sequence
of pebbling moves with certain properties. In many contexts, considering
a minimum sequence of pebbling moves with the properties in question
provides us with additional structure. For example, the no cycle lemma
implies that a minimum sequence of pebbling moves witnessing that
$p$ covers $q$ must be acyclic. 

We define a \emph{proper sink} to be a sink with indegree at least
one. 

\begin{lem}
[Minimum Signatures Lemma]\label{lem:minimum-signatures} Let $\sigma$
be a minimum sequence of pebbling moves in $G$ under $p$ which places
at least $k$ pebbles on $r$ with $p(r)\leq k$. If $D$ is the signature
of $\sigma$, then $D$ is acyclic, contains no proper sinks except
possibly $r$, the outdegree of $r$ is $0$, and the indegree of
$r$ is $k-p(r)$.
\end{lem}
\begin{proof}
By the no cycle lemma, $D$ is acyclic, or else we obtain a shorter
sequence of pebbling moves placing at least $k$ pebbles on $r$.
By \prettyref{lem:large-outdegree-many-pebbles}, the outdegree of
$r$ is zero, or again we obtain a shorter sequence. 

Because $d_{D}^{+}(r)=0$, we have $\balance(D,p,r)=p(r)+d_{D}^{-}(r)$.
Together with $\balance(D,p,r)\geq k$, we have that $d_{D}^{-}(r)\geq k-p(r)$.
If $d_{D}^{-}(r)>k-p(r)$, then $\balance(D,p,r)>k$. Obtain $D'$
from $D$ by deleting one edge into $r$. Notice that $D'$ satisfies
the balance condition and furthermore $\balance(D',p,r)=\balance(D,p,r)-1\geq k$.
By the acyclic orderability characterization, we obtain a shorter
sequence. 
\end{proof}
If we are interested in minimum sequences of pebbling moves that place
$k$ pebbles on some $r$ in a set $R$ of target vertices, the structure
of these sequences is further constrained. Not only do their signatures
obey the conditions found in the minimum signatures lemma, but the
outdegree of each vertex in $S$ is bounded.

\begin{lem}
\label{lem:minimum-signatures-R}Let $\sigma$ be a sequence of pebbling
moves in $G$ under $p$ that places at least $k>0$ pebbles on a
vertex $r\in R$ which, among all sequences placing at least $k$
pebbles on some vertex in $R$, minimizes the total number of pebbling
moves. Let $D$ be the signature of $\sigma$. For each $v\in R$,
we have that the outdegree of $v$ is less than $k/2$.
\end{lem}
\begin{proof}
Observe that $D$ is acyclic, or else we contradict the no cycle lemma.
Suppose for a contradiction that there is $v\in R$ with $d_{D}^{+}(v)\geq k/2$.
Because $k>0$, we have $d_{D}^{+}(v)>0$ and so \prettyref{lem:large-outdegree-many-pebbles}
yields a shorter sequence of pebbling moves placing at least $k$
pebbles on $v$, a contradiction.
\end{proof}

\section{Pebble Reachability}

Recall that the pebbling number of a graph $\pi(G)$ is the minimum
$k$ such that every vertex is reachable under every distribution
of size $k$. It is natural, then, to explore the decision problem
that results when we fix a particular distribution and target vertex;
that is, given $G$, $p$, and $r$, is $r$ reachable? We call this
problem $\PRlong$, or $\PR$ for short. As we show, $\PR$ is NP-complete,
even when the inputs are restricted so that $G$ is bipartite, has
maximum degree three, and each vertex starts with at most two pebbles. 

Analogously, fixing the distribution in the cover pebbling number
$\gamma(G)$ yields another decision problem: given $G$ and $p$,
does $p$ cover the unit distribution? We call this problem $\PClong$,
abbreviated $\PC$. Although deciding whether $\gamma(G)\leq k$ is
possible in polynomial time \cite{cover-thm-1,cover-thm-2}, $\PC$
is NP-complete. 

A sequence of pebbling moves $\sigma$ is \emph{nonrepetitive} if
for every (unordered) pair of vertices $\{ u,v\}$, $\sigma$ contains
at most one pebbling move between the vertices $u$ and $v$. Similarly
to $\PR$, we may ask, given $G$, $p$, and $r$, whether $r$ is
reachable via a nonrepetitive sequence of pebbling moves. We call
this language $\NPR$ (nonrepetitive pebble reachability). We show
that $\NPR$ is NP-complete. Our reduction is from a restricted form
of $\TSAT$ whose instances $\phi$ are all in a canonical form.

\begin{defn}
\label{def:Canonical-form}A $3\textrm{CNF}$ formula $\phi$ is in
\emph{canonical form} if 
\begin{enumerate}
\item $\phi$ has at least 2 clauses,
\item each clause contains 2 or 3 variables,
\item each variable appears at most 3 times in $\phi$,
\item each variable appears either once or twice in its positive form, and
\item each variable appears exactly once in its negative form
\end{enumerate}
\end{defn}
It is well known that $\TSAT$ remains NP-complete when (1-3) are
required. Suppose $\phi$ is a $\TSAT$ formula which satisfies (1-3)
but not necessarily (4) or (5). Indeed, if a variable $x$ always
appears in its positive (negative) form in $\phi$, we obtain a simpler,
equivalent formula by setting $x$ to true (false), thus removing
all clauses containing $x$ ($\overline{x}$). If $x$ appears twice
in its negative form, we simply switch all negative occurrences of
$x$ to positive occurrences and all positive occurrences of $x$
to negative occurrences. In this way, we obtain an equivalent formula
satisfying all of the above. We define $\RTSAT$ to be this restricted
form of $\TSAT$.

Our reduction from $\RTSAT$ to $\NPR$ employs several simple gadgets.
The AND gadget is a vertex $v$ that has two input edges and one output
edge; initially, $v$ is given zero pebbles. Notice that if $\sigma$
is nonrepetitive and contains a pebbling move from $v$ along the
output edge, then $\sigma$ must contain pebbling moves into $v$
along both input edges. The OR gadget is identical, except that $v$
is initially given a single pebble. In this case, if $\sigma$ is
nonrepetitive and contains a pebbling move from $v$ along the output
edge, then $\sigma$ must contain a pebbling move into $v$ along
one if the input edges. Using $2$-ary AND (OR) gadgets, one easily
constructs $k$-ary AND (OR) gadgets. 

The variable gadget is a path $v_{1}v_{2}v_{3}$ of length three.
The endpoint vertices $\left\{ v_{1},v_{3}\right\} $ are initially
given two pebbles, and the internal vertex $v_{2}$ is initially given
zero pebbles. The endpoint vertices correspond to the positive occurrence(s)
of the variable in $\phi$, and the internal vertex corresponds to
the negative occurrence of the variable in $\phi$. The variable gadget
has two or three output edges, depending upon how many times the corresponding
variable appears in $\phi$. If $x_{i}$ appears three times in $\phi$,
then its associated variable gadget $X_{i}$ has three output edges,
one incident to each $v_{i}$. If $x_{i}$ appears twice in $\phi$,
then $X_{i}$ has two output edges, one incident to each of $v_{1}$
and $v_{2}$. We say that the output edges incident to $v_{1}$ and
$v_{3}$ are \emph{positive output edges} and the output edge incident
to $v_{2}$ is the \emph{negative output edge}.

Given an instance $\phi$ of $\RTSAT$, we construct $G=G^{\NPR}(\phi)$
as follows. For each variable $x_{i}$ in $\phi$, we introduce a
variable gadget $X_{i}$ in $G$. For each clause $c_{j}$ containing
$k\in\left\{ 2,3\right\} $ variables, we introduce a $k$-ary OR
gadget $C_{j}$. The output edges of the $X_{i}$ are identified with
the input edges of the $C_{j}$ in the natural way: if $x_{i}$ appears
in $c_{j}$, a positive output edge of $X_{i}$ is identified with
an input edge of $C_{j}$, and if $\overline{x_{i}}$ appears in $c_{j}$,
the negative output edge of $X_{i}$ is identified with an input edge
of $C_{j}$. The output edges of the $C_{j}$ are connected to the
input edges of an $m$-ary AND gadget $A$, where $m$ is the number
of clauses in $\phi$. Finally, the output edge of $A$ is connected
to the target vertex $r$.

\begin{example*}
If $\phi=(w\vee x)\wedge(w\vee\overline{x})\wedge(\overline{w}\vee y\vee z)\wedge(x\vee\overline{y}\vee\overline{z})$,
then $G^{\NPR}(\phi)$ appears in \prettyref{fig:example-GNPR}.%
\begin{figure}
\begin{center}\includegraphics{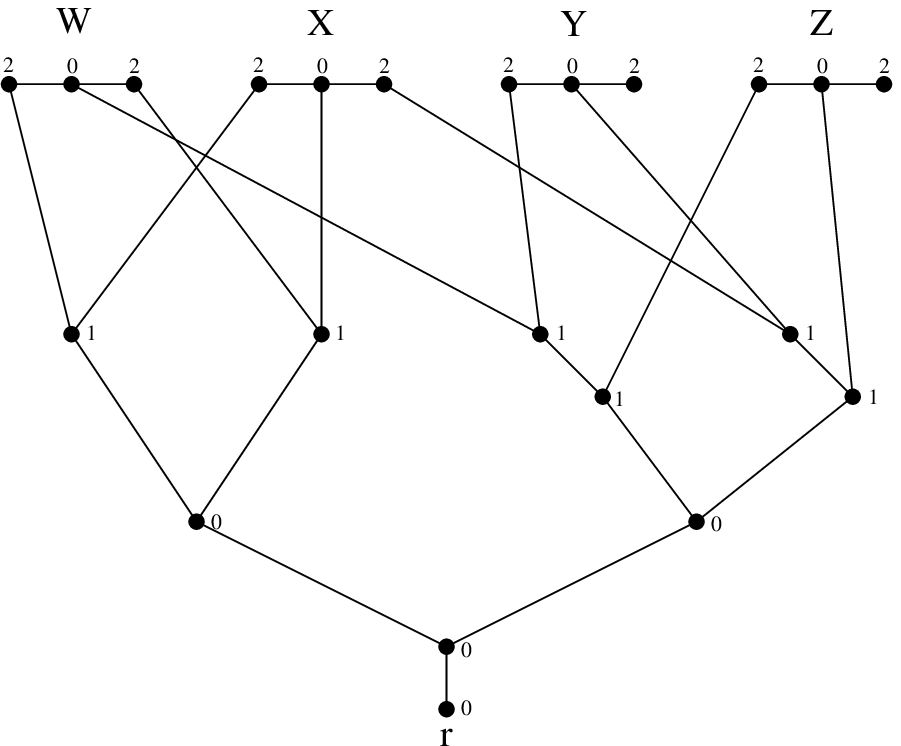}\end{center}

\caption{If $\phi=(w\vee x)\wedge(w\vee\overline{x})\wedge(\overline{w}\vee y\vee z)\wedge(x\vee\overline{y}\vee\overline{z})$,
then $G^{\NPR}(\phi)$ appears above. \label{fig:example-GNPR}}
\end{figure}

\end{example*}
\begin{prop}
\label{pro:GNPR-size}Let $\phi$ be an instance of $\RTSAT$ with
$n$ variables and $m$ clauses. Then $G^{\NPR}(\phi)$ has at most
$O(n+m)$ vertices. 
\end{prop}
\begin{thm}
\label{thm:NPR-is-NP-hard.} $\NPR$ is NP-complete, even when $G$
has maximum degree three and each vertex starts with at most two pebbles.
\end{thm}
\begin{proof}
It is immediate that $\NPR$ is in NP. Let $\phi$ be an instance
of $\RTSAT$ and let $G=G^{\NPR}(\phi)$. Observe that each vertex
in $G$ starts with at most two pebbles and the maximum degree in
$G$ is three. 

We claim that $\phi$ is satisfiable if and only if there is a nonrepetitive
sequence of pebbling moves which ends with a pebble on $r$. Suppose
that $\phi$ is satisfiable via $f\,:\,\left\{ x_{1},\ldots,x_{n}\right\} \to\left\{ \textrm{true},\textrm{false}\right\} $.
We construct a nonrepetitive sequence of pebbling moves which ends
with a pebble on $r$ as follows. For each variable $x_{i}$ with
$f(x_{i})=\textrm{false}$, we make a pebbling move from each endpoint
of $X_{i}$ to the interior vertex of $X_{i}$. Notice that after
executing these pebbling moves, for each $x_{i}$ with $f(x_{i})=\textrm{true}$,
we have two pebbles on each endpoint of $X_{i}$ and for each $x_{i}$
with $f(x_{i})=\textrm{false}$, we have two pebbles on the interior
vertex of $X_{i}$. Because $f$ satisfies $\phi$, each clause gadget
$C_{i}$ has some input edge which is incident to a vertex in a variable
gadget with two pebbles. By construction, each vertex in a variable
gadget is incident to at most one clause gadget input edge; therefore
we are able to make pebbling moves into each clause gadget $C_{i}$.
By the construction of our clause gadgets, we are then able to make
pebbling moves out of each clause gadget and, by construction, along
each of the inputs to the $m$-ary AND gadget. It follows that we
are able to make a pebbling move along the output of our AND gadget,
which places a pebble on $r$. It is easily observed that our sequence
of pebbling moves is nonrepetitive.

Conversely, suppose that $\sigma$ is nonrepetitive sequence of pebbling
moves which ends with a pebble on $r$. We construct a satisfying
assignment $f$ as follows. Because $\sigma$ contains a pebbling
move across the output of the AND gadget $A$, it follows that $\sigma$
contains pebbling moves across the output of each clause gadget $C_{i}$.
Hence, for each clause gadget $C_{i}$, $\sigma$ contains a pebbling
move across an input edge $e_{i}$ of $C_{i}$. If $e_{i}$ is incident
to an endpoint of $X_{j}$, then we set $f(x_{j})=\textrm{true}$;
otherwise, if $e_{i}$ is incident to the interior vertex of $X_{j}$,
we set $f(x_{j})=\textrm{false}$. We claim that we do not attempt
to set both $f(x_{j})=\textrm{true}$ and $f(x_{j})=\textrm{false}$.
Indeed, if we set $f(x_{j})=\textrm{false}$, then $\sigma$ contains
a pebbling move out of the interior vertex $v$ of $X_{j}$ along
an input edge to some clause gadget. Because $\sigma$ is nonrepetitive,
$v$ starts with zero pebbles, and $v$ has degree three, it must
be that $\sigma$ contains pebbling moves from each of the endpoints
in $X_{j}$ into $v$. Because each endpoint of $X_{j}$ starts with
only two pebbles and $\sigma$ is nonrepetitive, the moves into $v$
are the only pebbling moves which originate from the endpoints of
$X_{j}$. Therefore $\sigma$ does not contain a pebbling move out
of an endpoint of $X_{j}$ along an input edge of a clause gadget,
and hence we never attempt to set $f(x_{j})=\textrm{true}$. If the
truth values for any variables remain unset, we set them arbitrarily.
Now $f$ witnesses that $\phi$ is satisfiable.
\end{proof}
One of the major tools available to us when designing interesting
graph pebbling problems is the path; on a path, the pebbling moves
available to us are rather limited. If we are in a situation where
we need not concern ourselves with pebbling in cycles, then our options
on a path become even more limited. Furthermore, if the path is long,
it may be difficult to pebble across. Before using paths to reduce
$\NPR$ to $\PR$, we explore some basic properties.

\begin{lem}
\label{lem:path-edge-multiplicity}Let $G$ be a graph which contains
an induced path $P=v_{0}\ldots v_{n+1}$ containing $n+2$ vertices,
and suppose that each of the $n$ internal vertices in $P$ contains
$c$ pebbles. Let $D$ be an acyclic signature of a sequence of pebbling
moves so that the edge $v_{1}v_{0}$ has multiplicity $a_{0}\geq c$.
Then the multiplicity of $v_{n+1}v_{n}$ is at least $2^{n}(a_{0}-c)+c$.
\end{lem}
\begin{proof}
Observe that the claim is trivial if $a_{0}=0$; we assume that $a_{0}\geq1$.
For $1\leq i\leq n$, let $a_{i}$ be the multiplicity of $v_{i+1}v_{i}$.
We claim that for all $1\leq i\leq n$, we have that 
\begin{enumerate}
\item $a_{i}+c\geq2a_{i-1}$, and
\item $a_{i}\geq a_{0}$.
\end{enumerate}
Suppose for a contradiction that $i\geq1$ is the least integer for
which (1) or (2) fails, and consider the vertex $v_{i}$. By our selection
of $i$, $a_{i-1}\geq a_{0}$ and therefore $D$ requests at least
$a_{0}$ pebbling moves out of $v_{i}$ along edge $v_{i}v_{i-1}$.
Because $a_{0}\geq c$ and $a_{0}\geq1$, we have that $2a_{0}>c$;
hence, by the balance condition at $v_{i}$, the indegree of $v_{i}$
in $D$ is at least one. Because $D$ is acyclic, $D$ contains no
edges of the form $v_{i-1}v_{i}$. Because $v_{i}$ is an internal
vertex in an induced path in $G$, the only other edge incident to
$v_{i}$ is $v_{i}v_{i+1}$. It follows that the indegree of $v_{i}$
in $D$ is exactly the multiplicity of $v_{i+1}v_{i}$, and so the
indegree of $v_{i}$ in $D$ is $a_{i}$. Therefore the balance condition
at $v_{i}$ implies that $a_{i}+c\geq2a_{i-1}$, which together with
$a_{0}\geq c$ and $a_{i-1}\geq a_{0}$, implies $a_{i}\geq a_{0}$. 

Solving our recurrence in (1), we find that $a_{i}\geq2^{i}(a_{0}-c)+c$. 
\end{proof}
We use our path lemma to argue that if we can pebble across a long
path several times, than we can place many pebbles on the originating
endpoint of the path. Together with \prettyref{lem:large-outdegree-many-pebbles},
we obtain the following corollary.

\begin{cor}
\label{cor:many-pebbles-begining-of-path}Under the assumptions of
\prettyref{lem:path-edge-multiplicity}, there exists $D'\subseteq D$
such that $D'$ is orderable and $\balance(D',p,v_{n+1})\geq2^{n+1}(a_{0}-c)+2c$.
If in addition we have $d_{D}^{+}(v_{n+1})>0$, then we may take $D'$
to be a proper subgraph of $D$.
\end{cor}
Our reduction used the notion of nonrepetitive sequences of pebbling
moves. In fact, there is a natural correspondence between the nonrepetitive
sequences of pebbling moves in a graph $G$ and (arbitrary) sequences
of pebbling moves in another graph $\mathcal{S}(G,\alpha)$. 

\begin{defn}
\label{def:subdivide-graph}We obtain $\mathcal{S}(G,\alpha)$ from
$G$ by replacing each edge in $G$ with a path containing $\alpha$
internal vertices, so that $d_{\mathcal{S}(G,\alpha)}(u,v)=\left(1+\alpha\right)d_{G}(u,v)$
for any $u,v$ in $G$. We call these paths \emph{one use paths}. 
\end{defn}
As our next lemma shows, the correspondence holds whenever $\alpha$
is sufficiently large with respect to the number of pebbles in $G$. 

\begin{lem}
\label{lem:one-use-edges}Fix a graph $G$ and a parameter $t\geq0$.
Suppose that $\alpha\geq\max\left\{ \lg2t,\,4\lg e(G)\right\} $ and
let $H=\mathcal{S}(G,\alpha)$. Let $p$ be a pebble distribution
on $G$ of size at most $t$ and define a pebble distribution $q$
on $H$ so that $q$ and $p$ agree on $V(G)$ and $q$ assigns one
pebble each to the internal vertices of $H$'s one use paths. We have
the following claims.
\begin{enumerate}
\item If $\sigma$ is a nonrepetitive sequence of pebbling moves in $G$,
then there exists a sequence of pebbling moves $\sigma'$ in $H$
such that $p_{\sigma}$ and $q_{\sigma'}$ agree on $V(G)$. 
\item Conversely, if $\sigma$ is a sequence of pebbling moves in $H$,
then there exists a nonrepetitive sequence of pebbling moves $\sigma'$
in $G$ such that $p_{\sigma'}(v)\geq q_{\sigma}(v)$ for all $v$
in $G$.
\end{enumerate}
\end{lem}
\begin{proof}
Claim 1 is clear. Suppose that $\sigma$ is a sequence of pebbling
moves in $H$. By the no cycle lemma, we may assume without loss of
generality that the signature $D$ of $\sigma$ is acyclic. We define
a digraph $D'$ with vertex set $V(G)$ as follows. Let $uv$ be an
edge in $G$ and let $u=w_{0}\ldots w_{\alpha+1}=v$ be the corresponding
one use path in $H$. The multiplicity of the edge $uv$ in $D'$
is the multiplicity of the edge $w_{\alpha}w_{\alpha+1}$ in $D$.
Because $D$ is acyclic, the balance condition implies that if $D$
contains the edge $w_{\alpha}w_{\alpha+1}$, then $D$ contains all
edges $w_{k}w_{k+1}$. It follows that $D'$ is also acyclic. It is
easily seen that $\balance(D',p,v)\geq\balance(D,q,v)$ for each $v$
in $D'$. By the acyclic orderability characterization, we obtain
a sequence of pebbling moves $\sigma'$ such that $p_{\sigma'}(v)\geq q_{\sigma}(v)$
for all $v$ in $G$. It remains to show that $D'$ has no edges of
multiplicity at least two, so that $\sigma'$ is necessarily nonrepetitive. 

Suppose for a contradiction that $uv$ is an edge in $D'$ with multiplicity
at least two; again, let $u=w_{0}\ldots w_{\alpha+1}=v$ be the corresponding
one use path in $H$. It follows that $w_{\alpha}w_{\alpha+1}$ has
multiplicity at least two in $D$. Recalling that $q$ assigns each
of the internal vertices $w_{i}$ one pebble, \prettyref{lem:path-edge-multiplicity}
implies that the multiplicity of $w_{0}w_{1}$ is at least $2^{\alpha}+1$.
Because each pebbling move reduces the total number of pebbles by
one, certainly the size of $q$ is at least $2^{\alpha}+2$. But $\left|q\right|=\left|p\right|+\alpha e(G)$
and together with $t\leq2^{\alpha-1}$ and $\alpha e(G)\leq2^{\alpha-1}$,
we obtain a contradiction. 
\end{proof}
\begin{cor}
$\PRlong$ is NP-complete, even when $G$ is bipartite, has maximum
degree three, and each vertex starts with at most two pebbles.
\end{cor}
\begin{proof}
By the no cycle lemma and the acyclic orderability characterization,
$\PR$ is in NP. We reduce $\NPR$ to $\PR$ as follows. Consider
a graph $G$ with maximum degree three, a distribution of pebbles
$p$ which places at most two pebbles on each vertex in $G$, and
a target vertex $r$. Let $\alpha$ be the least odd number larger
than $\max\left\{ \lg2\left|p\right|,\,4\lg e(G)\right\} $. Our reduction
outputs $H=\mathcal{S}(G,\alpha)$ with pebble distribution $q$ as
in \prettyref{lem:one-use-edges} and target vertex $r$. Observe
that $H$ is bipartite, has maximum degree three, and each vertex
starts with at most two pebbles. By \prettyref{lem:one-use-edges},
$r$ is reachable via a nonrepetitive sequence of pebbling moves in
$G$ if and only if $r$ is reachable in $H$. 
\end{proof}
Let $\phi$ be an instance of $\RTSAT$. We define $G^{\PR}(\phi)=\mathcal{S}(G^{\NPR}(\phi),\alpha)$
with $\alpha$ chosen as in our corollary; that is, $G^{\PR}$ is
the composition of our reduction from $\RTSAT$ to $\NPR$ and our
reduction from $\NPR$ to $\PR$. 

\begin{cor}
$\PClong$ is NP-complete, even when $G$ is bipartite, has maximum
degree three, and each vertex starts with at most three pebbles.
\end{cor}
\begin{proof}
By the no cycle lemma and the acyclic orderability characterization,
we have that $\PC$ is in NP. We reduce $\PR$ to $\PC$ as follows.
Let $G$ be a graph with pebble distribution $p$ and target vertex
$r$. Define a new distribution $q$ of pebbles so that $q(v)=p(v)+1$
for all $v\neq r$ and $q(r)=p(r)$. We claim that $r$ is reachable
under $p$ if and only if $q$ covers the unit distribution. The forward
direction is clear. 

Suppose that $\sigma$ is a minimum sequence of pebbling moves witnessing
that $q$ covers the unit distribution, and let $D$ be the signature
of $\sigma$. By the no cycle lemma, $D$ is acyclic. Because $\balance(D,q,v)\geq1$,
we have that $\balance(D,p,v)\geq0$ for all $v$ and $\balance(D,p,r)\geq1$.
It follows from the acyclic orderability characterization that $D$
is orderable under $p$. Together with $\balance(D,p,r)\geq1$, we
have that $r$ is reachable under $p$.
\end{proof}
As we have seen, $\PRlong$ is NP-complete, even under some restrictions
of the inputs. However, as we now observe, if we restrict $G$ to
be a tree, then we can solve $\PRlong$ in polynomial time using a
simple greedy strategy. A \emph{greedy pebbling move} is a pebbling
move $uv$ such that $d(v,r)<d(u,r)$. The \emph{greedy pebbling strategy}
arbitrarily makes greedy pebbling moves until no greedy pebbling move
is possible. 

\begin{prop}
[Greedy Tree Lemma]\label{pro:tree-greedy-strategy}In a tree $T$
with target $r$, the maximum number of pebbles that can be placed
on $r$ is achieved with the greedy pebbling strategy.
\end{prop}
\begin{proof}
Suppose for a contradiction that under $p$, it is possible to place
$k$ pebbles on $r$, but if we make the greedy pebbling move $uv$,
it is no longer possible to place at least $k$ pebbles on $r$. Let
$\sigma$ be a minimum sequence of pebbling moves placing $k$ pebbles
on $r$, and let $D$ be the signature of $\sigma$. By the no cycle
lemma, $D$ is acyclic. If $D$ contains the edge $uv$, then the
acyclic orderability characterization implies that $D-uv$ is orderable
under $p_{uv}$, implying it is possible to place $k$ pebbles on
$r$ even after pebbling $uv$. Otherwise, if $D$ does not contain
the edge $uv$, then $d^{+}(u)=0$, or else $D$ contains a proper
sink other than $r$, contradicting the minimal signatures lemma.
Therefore $\sigma$ does not contain any pebbling moves out of $u$,
and so $uv$ followed by $\sigma$ is a legal sequence of pebbling
moves placing at least $k$ pebbles on $r$. 
\end{proof}

\section{Complexity of Optimal Pebbling Number}

Recall that the optimal pebbling number $\piopt(G)$ of a graph $G$
is the least number $k$ such that every vertex is reachable under
some distribution of size $k$. We define $\OPNlong$ (abbreviated
$\OPN$) to be the problem of deciding, given $G$ and $k$, whether
$\piopt(G)\leq k$. In this section, we show that $\OPN$ is NP-complete.
We observe that $\OPN$ is in NP; indeed, we may witness that $\piopt(G)\leq k$
by providing a distribution $p$ of size $k$ and, for each $r$,
the signature $D_{r}$ of a sequence of pebbling moves showing that
$r$ is reachable. More care is needed to establish that $\OPN$ is
NP-hard. As in our proof that $\PR$ is NP-hard, we establish that
$\OPN$ is NP-hard through an intermediate decision problem. 

Let $G$ be a graph and $p$ be a distribution of pebbles to $G$.
A vertex $r$ is \emph{determinative} if $r$ is reachable under $p$
implies that every vertex in $G$ is reachable under $p$. Informally,
if $r$ is determinative, then no vertices in $G$ are more difficult
to pebble than $r$. Our intermediate decision problem is $\PRlong$
with the added restriction that $r$ is determinative. We call this
problem $\DPR$ (determinative pebble reachability). 

\begin{prop}
\label{pro:-DPR-is-NP-hard}$\DPR$ is NP-complete, even when each
vertex starts with at most two pebbles.
\end{prop}
\begin{proof}
Because $\PRlong$ is in NP, it is immediate that $\DPR$ is in NP
as well. We show that our reduction $G^{\PR}$ from $\RTSAT$ to $\PR$
actually produces an instance of $\DPR$. Let $\phi$ be an instance
of $\RTSAT$, and let $G=G^{\PR}(\phi)$ with distribution $p$ and
target $r$. We show that $r$ is determinative. Suppose that it is
possible to place a pebble on $r$, or equivalently that $\phi$ is
satisfiable. Consider a vertex $v\in G$. If $v$ is an internal vertex
in a one use path introduced in our reduction from $\NPR$ to $\PR$,
then $v$ begins with one pebble and so $v$ is reachable trivially. 

It remains to consider the case that $v$ is a vertex introduced in
our reduction from $\RTSAT$ to $\NPR$, so that $v$ is either in
an OR gadget, a variable gadget, an AND gadget, or $v=r$. If $v$
is in an OR gadget, then $v$ begins with a pebble. If $v$ is an
endpoint of a variable gadget, then $v$ begins with two pebbles.
If $v$ is the interior vertex of a variable gadget, then we may place
a pebble on $v$ by pebbling from either of the endpoints (which start
with two pebbles) across the one use path. Otherwise, if $v$ is in
an AND gadget or $v=r$, then we use the satisfiability of $\phi$
to place a pebble on $v$.
\end{proof}
Before we are able to present our reduction from $\DPR$ to $\OPN$,
we require some technical lemmas. The following weighting argument
is well known and is a fundamental tool in graph pebbling.

\begin{prop}
[Standard Weight Equation]\label{pro:weight-eqn}Let $G$ be a graph
with distribution $p$ and target vertex $r$, and let $a_{i}$ be
the number of pebbles at distance $i$ from $r$. If it is possible
to place $s$ pebbles on $r$, then we have $\sum_{i\geq0}2^{-i}a_{i}\geq s$.
\end{prop}
\begin{proof}
Observe that it is not possible to make a pebbling move which increases
the sum $\sum_{i\geq0}2^{-i}a_{i}$. 
\end{proof}
The following graph will be useful to us in two different contexts:
first, as a gadget, and secondly in establishing the correctness of
our reduction from $\DPR$ to $\OPN$.

\begin{defn}
We define $\starg(\alpha,\beta)$ to be the result of replacing each
edge in $K_{1,\beta}$ with a path of length $\alpha$, so that $\starg(\alpha,\beta)$
has $\alpha\beta$ edges. Equivalently, $\starg(\alpha,\beta)=\mathcal{S}(K_{1,\beta},\alpha-1)$. 
\end{defn}
\begin{example*}
$\starg(3,5)$ appears in \prettyref{fig:starg-example}. %
\begin{figure}
\begin{center}\includegraphics[%
  scale=0.4]{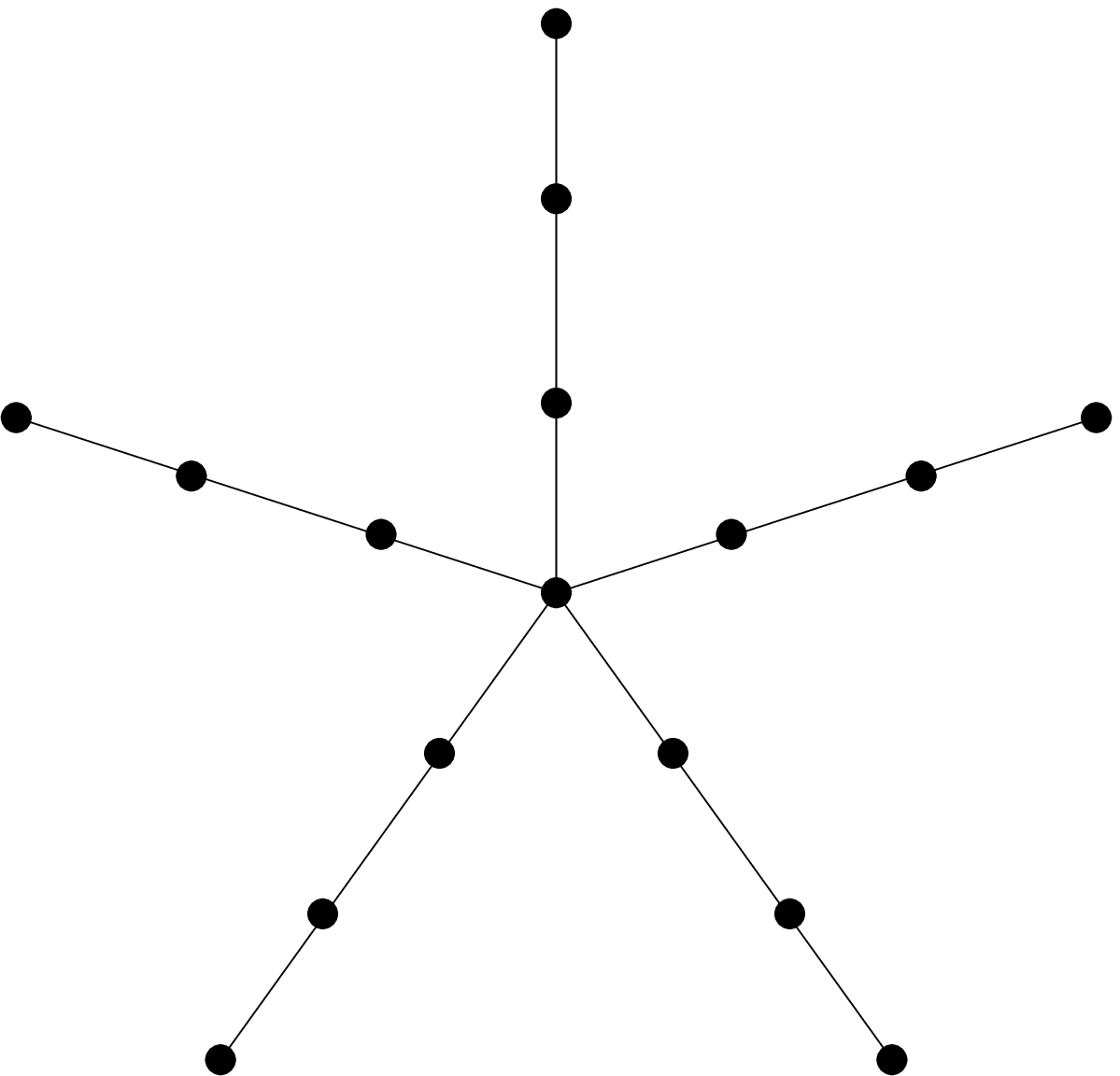}\end{center}

\caption{$\starg(3,5)$\label{fig:starg-example}}
\end{figure}

\end{example*}
Our reduction from $\DPR$ to $\OPN$ produces a graph whose global
structure is similar to that of $\starg(\cdot)$. Our instance of
$\DPR$ plays the role of the center vertex, and the gadgets that
we add play the role of the leaves. When we argue the correctness
of our reduction, we apply the following lemma to limit the pebble
distributions that we must consider. The lemma shows how, despite
the simplicity of the standard weight equation, it yields nontrivial
results.

\begin{lem}
\label{lem:star-lemma} Fix $\alpha\geq1$ and $\beta\geq2$. Let
$p$ be a distribution of $\beta2^{\alpha}$ pebbles to $\starg(\alpha,\beta)$
with the property that for each leaf $l$ in $\starg(\alpha,\beta)$,
it is possible to place $2^{\alpha}$ pebbles on $l$. If $2(\beta^{2}+1)<2^{\alpha}$,
then $p$ is the distribution which places $2^{\alpha}$ pebbles on
each leaf and zero pebbles on the other vertices.
\end{lem}
\begin{proof}
Let $v$ be the center vertex of $\starg(\alpha,\beta)$, and for
each $0\leq i\leq\alpha$, let $a_{i}$ be the number of pebbles at
distance $i$ from $v$. For each $l$, it is possible to place $2^{\alpha}$
pebbles on $l$ and \prettyref{pro:weight-eqn} yields an equation;
we sum these equations. Because there are $\beta$ leaves, we obtain
$\beta2^{\alpha}$ on the right hand side. A pebble at distance $i$
from $v$ is at distance $\alpha-i$ from its closest leaf and $\alpha+i$
from all other leaves. It follows that pebbles at distance $i$ from
$v$ contribute $1/2^{\alpha-i}+(\beta-1)/2^{\alpha+i}$ to the left
hand side of the equation. We obtain \[
\sum_{i=0}^{\alpha}\left(\frac{1}{2^{\alpha-i}}+\frac{\beta-1}{2^{\alpha+i}}\right)a_{i}\geq\beta2^{\alpha}\]
 and after some simplification, \[
\sum_{i=0}^{\alpha}\left(2^{i}+\frac{\beta-1}{2^{i}}\right)a_{i}\geq\beta4^{\alpha}.\]
 Let $f(x)=2^{x}+(\beta-1)2^{-x}$, so that pebbles at distance $i$
contribute $f(i)$ to the left hand side. Analyzing the derivative
$f'(x)=\ln2\left(2^{x}-(\beta-1)2^{-x}\right)$, we find that $f'(x)=0$
has one solution, namely $x_{0}=\log_{4}(\beta-1)$. Furthermore,
for $x>x_{0}$, we have $f'(x)>0$ and for $x<x_{0}$, we have $f'(x)<0$.
It follows that $f(x)$ has a global minimum at $x=x_{0}$, $f(x)$
is decreasing on $(-\infty,x_{0}]$, and $f(x)$ is increasing on
$[x_{0},\infty)$. 

Let $m=\sum_{i=0}^{\alpha-1}a_{i}$ be the number of pebbles not at
distance $\alpha$ from $v$; we show that $m<1$, implying that $m=0$.
Noting that $a_{\alpha}=\beta2^{\alpha}-m$, we have that \[
\left(\max_{0\leq i\leq\alpha-1}f(i)\right)m+f(\alpha)\left(\beta2^{\alpha}-m\right)\geq\beta4^{\alpha}.\]
Because of the monotonicity properties of $f$, we have $\max_{i=0}^{\alpha-1}f(i)\in\left\{ f(0),f(\alpha-1)\right\} $.
Because $2(\beta^{2}+1)<2^{\alpha}$, certainly $2\beta<2^{\alpha}$
and therefore \[
f(0)=\beta<2^{\alpha-1}\leq2^{\alpha-1}+(\beta-1)2^{1-\alpha}=f(\alpha-1)\]
 It follows that $\max_{0\leq i\leq\alpha-1}f(i)=f(\alpha-1)$, and
after substitution and further simplification, we obtain \[
m\leq\frac{\beta2^{\alpha}\left(f(\alpha)-2^{\alpha}\right)}{f(\alpha)-f(\alpha-1)}.\]
 Substituting our formula for $f(\alpha)$ into the numerator yields
\[
m\leq\frac{\beta(\beta-1)}{f(\alpha)-f(\alpha-1)}\leq\frac{\beta^{2}}{f(\alpha)-f(\alpha-1)}.\]
 Recall that $2(\beta^{2}+1)<2^{\alpha}$, which implies that $\beta^{2}<2^{\alpha-1}-1$.
Observe that $f(\alpha)-f(\alpha-1)=2^{\alpha-1}-(\beta-1)/2^{\alpha}$.
Because $\beta-1<2(\beta^{2}+1)<2^{\alpha}$, we have that $f(\alpha)-f(\alpha-1)>2^{\alpha-1}-1$.
It follows that $\beta^{2}<f(\alpha)-f(\alpha-1)$, implying that
$m<1$ as required. 

It follows that $p$ places every pebble at distance $\alpha$ from
$v$. It remains to show that $p$ places $2^{\alpha}$ pebbles on
each leaf. Fix an arbitrary leaf $l$, and let $n$ be the number
of pebbles that $p$ places on $l$. Applying the standard weight
equation to $l$, we have that \[
n+\frac{\beta2^{\alpha}-n}{2^{2\alpha}}\geq2^{\alpha}.\]
 After simplification, we obtain that \begin{eqnarray*}
n & \geq & 2^{\alpha}-\frac{2^{\alpha}(\beta-1)}{4^{\alpha}-1}.\end{eqnarray*}
 Similarly to the previous paragraph, we show that $n>2^{\alpha}-1$.
We have that \[
\frac{2^{\alpha}(\beta-1)}{4^{\alpha}-1}\leq\frac{2^{\alpha}(\beta-1)}{4^{\alpha}-2^{\alpha}}=\frac{\beta-1}{2^{\alpha}-1}.\]
 Because $\beta\leq2(\beta^{2}+1)<2^{\alpha}$, we have that $(\beta-1)/(2^{\alpha}-1)<1$
and hence $n>2^{\alpha}-1$ as required. Therefore $p$ assigns each
leaf at least $2^{\alpha}$ pebbles and the lemma follows.
\end{proof}
We now have the tools necessary to present our reduction from $\DPR$
to $\OPN$. Let $G$ be a graph with pebble distribution $p$ and
determinative target vertex $r$. Let $m=\left|p\right|$, let $\alpha=\left\lceil \lg\left(2(m^{2}+1)\right)\right\rceil $,
and let $\beta=2^{\alpha}m+2$. We construct a graph $H$ with the
property that $\piopt(H)\leq m2^{\alpha}$ if and only if $r$ is
reachable in $G$. 

We construct $H$ from $G$ by attaching a copy of $\starg(\alpha,\beta)$
to each pebble in $G$. That is, for each pebble on a vertex $u$,
we introduce a copy of $\starg(\alpha,\beta)$ and attach it to $u$
by identifying $u$ with one of the leaves of our copy of $\starg(\alpha,\beta)$.

\begin{lem}
$r$ is reachable in $G$ under $p$ if and only if $\piopt(H)\leq m2^{\alpha}$.
\end{lem}
\begin{proof}
($\implies$). Suppose $r$ is reachable. Define a distribution $q$
of $m2^{\alpha}$ pebbles to $H$ by placing $2^{\alpha}$ pebbles
at the centers of each of the $m$ copies of $\starg(\alpha,\beta)$
in $H$. Consider a vertex $v$ in $H$. If $v$ belongs to a copy
$S$ of $\starg(\alpha,\beta)$, then $v$ is at a distance of at
most $\alpha$ from the center of $S$; because the center of $S$
begins with $2^{\alpha}$ pebbles, $v$ is reachable. Otherwise, $v$
must be a vertex in $G$. Because $r$ is reachable and determinative
under $p$, to show that $v$ is reachable, it suffices to show that
$q$ covers $p$. But each star can contribute one pebble to the vertex
it shares with $G$, and so $q$ covers $p$. 

($\Longleftarrow$). Let $q$ be a distribution of $m2^{\alpha}$
pebbles to $H$ witnessing that $\piopt(H)\leq m2^{\alpha}$. We claim
that if $u$ is the center vertex of a copy $S$ of $\starg(\alpha,\beta)$,
then it is possible to place $2^{\alpha}$ pebbles on $u$ starting
from $q$. Indeed, because $S$ contains $\beta-1>m2^{\alpha}$ pendant
paths with endpoint $u$, there is some path to which $q$ assigns
no pebbles (except possibly at $u$). Let $w_{0}w_{1}\ldots w_{\alpha}$
be one such path with $w_{0}=u$. Because every vertex is reachable
under $q$, certainly $w_{\alpha}$ is reachable; let $D$ be a signature
of a minimum sequence of pebbling moves that places a pebble on $w_{\alpha}$.
Because $w_{\alpha}$ is a leaf and $q$ assigns no pebbles to $w_{\alpha}$,
$w_{\alpha-1}w_{\alpha}$ is an edge in $D$; therefore \prettyref{cor:many-pebbles-begining-of-path}
implies that we can place $2^{\alpha}$ pebbles on $u$. 

When a graph has a pebble distribution, contracting a set of vertices
$S$ changes the pebble distribution in the natural way: pebbles on
vertices in $S$ are collected at the vertex of contraction. Construct
$H'$ and pebbling distribution $q'$ from $H$ and $q$ by iteratively
applying the following contractions: 
\begin{enumerate}
\item Contract all vertices in $H$ that are also in $G$ to a single vertex
$v$
\item For each copy $S$ of $\starg(\alpha,\beta)$, contract the vertices
in $S$ that are at distance at least $\alpha$ from $v$
\end{enumerate}
Observe that $H'$ is exactly $\starg(\alpha,m)$, with center vertex
$v$. Because the contraction operation cannot make pebbling more
difficult, it is possible to place $2^{\alpha}$ pebbles on each leaf
in $H'$ starting from $q'$. Because $2(m^{2}+1)\leq2^{\alpha}$,
applying \prettyref{lem:star-lemma} to $H'=\starg(\alpha,m)$ implies
that $q'$ must assign $2^{\alpha}$ pebbles to each leaf of $H'$.
It follows that $q$ assigns $2^{\alpha}$ pebbles to each copy of
$\starg(\alpha,\beta)$ in $H$ in such a way that each pebble is
at distance at least $\alpha$ away from the vertices in $G$.

Let $E$ be the signature of a minimum sequence of pebbling moves
in $H$ starting from $q$ which places a pebble on $r$. Consider
a copy $S$ of $\starg(\alpha,\beta)$ attached to a vertex $u$ in
$G$. We claim that $E$ contains at most one edge from $S$ into
$u$. Indeed, if this were otherwise, then by \prettyref{cor:many-pebbles-begining-of-path}
there exists $E'\subseteq E$ which places at least $2\cdot2^{\alpha}$
pebbles on the center vertex of $S$. However, this is impossible
because $E$ is acyclic with edges from $S$ into $G$ and $q$ assigns
only $2^{\alpha}$ pebbles to $S$. 

Obtain $E'$ from $E$ by deleting all edges except those in $G$.
Because each vertex $u$ in $G$ receives a pebble from $p$ for every
attached copy of $\starg(\alpha,\beta)$, we have that $\balance(E',p,u)\geq\balance(E,q,u)$.
It follows from the acyclic orderability characterization that $E'$
is orderable under $p$; together with $\balance(E',p,u)\geq\balance(E,q,u)\geq1$,
we have that $r$ is reachable under $p$. 
\end{proof}
We conclude with this section's main theorem. 

\begin{thm}
$\OPNlong$ is NP-complete.
\end{thm}
\begin{proof}
We have already observed that $\OPN$ is in NP and exhibited a reduction
from $\DPR$ to $\OPN$. It remains to check that the size of $H$
and $m2^{\alpha}$ are not too large, so that our reduction is computable
in polynomial time. Let $n$ be the number of vertices in $G$. By
\prettyref{pro:-DPR-is-NP-hard}, we assume that $\left|p\right|=m$
is at most $2n$. Our reduction uses gadgets $\starg(\alpha,\beta)$
with $\alpha\leq\left\lceil \lg\left(2(4n^{2}+1)\right)\right\rceil $
and $\beta\leq2^{\alpha}m+2=O(n^{3})$. It follows that each gadget
$\starg(\alpha,\beta)$ has at most $O(n^{3}\log n)$ vertices. Because
we use at most $2n$ gadgets, $H$ contains a total of at most $n+2nO(n^{3}\log n)=O(n^{4}\log n)$
vertices. 
\end{proof}

\section{Complexity of Pebbling Number}

Although the optimal pebbling number has received some study, combinatorialists
have focused more attention on the pebbling number. Recall that the
$r$-pebbling number $\pi(G,r)$ is the minimum $k$ such that $r$
is reachable under every distribution of size $k$. Similarly, the
pebbling number $\pi(G)$ is the minimum $k$ such that every vertex
is reachable under every distribution of size $k$. It is clear from
the definitions that if $n$ is the number of vertices in $G$, then
$\piopt(G)\leq n\leq\pi(G)$. At first glance, it may not be clear
that $\pi(G)$ is well defined. In fact, if $G$ is not connected,
then we can place arbitrarily many pebbles in a single component and
we will not be able to place pebbles on vertices outside the component.
However, for connected graphs, $\pi(G)$ is well defined; we implicitly
assume that $G$ is connected. Indeed, if $d$ is the diameter of
$G$, every vertex is reachable provided that our distribution is
forced to place at least $2^{d}$ pebbles on some vertex. We record
this observation as a proposition.

\begin{prop}
\label{pro:PN-trivial-bound} Let $G$ be a graph with diameter $d$.
Then $\pi(G)\leq\left(2^{d}-1\right)n+1$. 
\end{prop}
We call the problem of deciding whether $\pi(G,r)\leq k$ $\RPNlong$
(abbreviated $\RPN$); similarly, we define $\PNlong$ (abbreviated
$\PN$) to be the problem of deciding whether $\pi(G)\leq k$. In
this section, we establish that $\PN$ and $\RPN$ are $\ptc$-complete.
First, note that both languages are in $\ptc$. Indeed, to decide
if $\pi(G)\leq k$, our machine need only check that for all distributions
$p$ of size $k$ and all target vertices $r$, there exists an orderable
digraph $D_{p,r}$ that places a pebble on $r$. The distributions
of size $k$, the target vertices, and the digraphs $D_{p,r}$ are
all describable using $\poly(n,\log k)$ bits. Further, $\POlong$
is in P. It follows that $\PN$ is in $\ptc$. A similar argument
shows that $\RPN$ is in $\ptc$.

The seminal $\ptc$-complete problem is a quantified version of $\TSAT$
whose instances consist of a $3\textrm{CNF}$ formula $\phi$ over
a set of universally quantified variables and a set of existentially
quantified variables (see \cite{complexity}). We say that $\phi$
is \emph{valid} if for every setting of the universally quantified
variables, there is a setting of the existentially quantified variables
which satisfies $\phi$. The decision problem $\AETSAT$ is to determine
whether $\phi$ is valid. 

Just as $\TSAT$ remains NP-complete when $\phi$ is restricted to
be in canonical form (recall \prettyref{def:Canonical-form}), $\AETSAT$
remains $\ptc$-complete when $\phi$ is restricted to be in canonical
form. We call this restriction $\RAETSAT$. 

We show that $\RPN$ is $\ptc$-complete by a reduction from $\RAETSAT$.
Whereas our reduction to $\OPN$ produces graphs $H$ with the property
that only one distribution can possibly succeed in witnessing $\piopt(H)\leq k$,
our reduction to $\RPN$ produces graphs with the property that almost
all distributions succeed in being able to place a pebble on $r$.
It is the rare {}``difficult'' distributions -- those which may
not allow a pebble to be placed on $r$ -- that correspond to settings
of the universally quantified variables in our $\RAETSAT$ formula.
Given a distribution of $k$ pebbles to the graph we produce, either
$r$ is easily reachable, or the distribution corresponds to a setting
$f$ of the universally quantified variables in $\phi$ and $r$ is
reachable if and only if $\phi$ is satisfiable under $f$. 

Our reduction from $\RAETSAT$ to $\RPN$ involves the construction
of several graphs, each building on the previous construction. We
refer to the $i$th graph we produce as $G_{i}=G_{i}(\phi)$. We present
the reduction with respect to a fixed instance $\phi$ of $\RAETSAT$.

\subsection{The Underlying Graph}

We obtain $G_{1}$ from $\phi$ by modifying $G^{\NPR}(\phi)$ slightly.
That is, for each universally quantified variable $x_{i}$ in $\phi$,
we remove both edges from the variable gadget $X_{i}$ in $G^{\NPR}(\phi)$
associated with $x_{i}$ and remove one pebble each from the endpoints
of $X_{i}$, so that the endpoints of $X_{i}$ start with one pebble
instead of two. (We leave intact variable gadgets $X_{j}$ corresponding
to existentially quantified variables $x_{j}$ in $\phi$.) Let $n_{1}=n(G_{1})$
be the number of vertices in $G_{1}$, let $e_{1}=e(G_{1})$, and
let $p_{1}$ be the distribution on $G_{1}$. The following definition
gives the correspondence between settings of the universally quantified
variables in $\phi$ and distributions of pebbles in $G_{1}$.

\begin{defn}
For each setting $f$ of the universally quantified variables in $\phi$,
let $p_{1,f}$ be the distribution of pebbles to $G_{1}$ given by
adding the following pebbles to $p_{1}$. For each $x_{i}$ with $f(x_{i})=\textrm{true}$,
add one pebble to each of the two vertices associated with positive
occurrences of $x_{i}$ in $\phi$. For each $x_{i}$ with $f(x_{i})=\textrm{false}$,
add two pebbles to the vertex associated with the negative instance
of $x_{i}$. 
\end{defn}
Observe that under any $p_{1,f}$, each vertex in $G_{1}$ contains
at most two pebbles. Our interest in $G_{1}$ under the distributions
$p_{1,f}$ is based on the following proposition, whose proof is similar
to that of \prettyref{thm:NPR-is-NP-hard.}. 

\begin{prop}
\label{pro:PN-corresp-G1} There is a nonrepetitive sequence of pebbling
moves which places a pebble on $r$ in $G_{1}$ starting from $p_{1,f}$
if and only if there is a setting of the existentially quantified
variables in $\phi$ which, together with $f$, satisfies $\phi$.
\end{prop}
Let $t$ be the number of pebbles in the $p_{1,f}$. Because $p_{1,f}$
assigns at most two pebbles to each vertex in $G_{1}$, $t\leq2n_{1}$.
We obtain $G_{2}$ from $G_{1}$ by setting $\alpha=\max\left\{ \lg2t,\,4\lg e_{1}\right\} $
and replacing each edge in $G_{1}$ with a path of length $\alpha+1$;
that is, $G_{2}=\mathcal{S}(G_{1},\alpha)$ (recall \prettyref{def:subdivide-graph}).
Let $n_{2}$ be the number of vertices in $G_{2}$.

Let $p_{2}$ be the distribution of pebbles to $G_{2}$ so that $p_{2}$
and $p_{1}$ agree on all vertices in $G_{1}$ and $p_{2}(v)=1$ for
all vertices $v$ introduced in our construction of $G_{2}$ from
$G_{1}$. Similarly, let $p_{2,f}$ be the distribution of pebbles
to $G_{2}$ so that $p_{2,f}$ and $p_{1,f}$ agree on all vertices
in $G_{1}$ and $p_{2,f}(v)=1$ for all vertices $v$ introduced in
our construction of $G_{2}$ from $G_{1}$. 

We call $G_{2}$ the \emph{underlying graph} and a distribution $p_{2,f}$
an \emph{underlying distribution}. Observe that by \prettyref{lem:one-use-edges},
there is a nonrepetitive sequence of pebbling moves which places a
pebble on $r$ in $G_{1}$ under $p_{1,f}$ if and only if there is
an arbitrary sequence of pebbling moves in $G_{2}$ under $p_{2,f}$
which places a pebble on $r$. Together with \prettyref{pro:PN-corresp-G1},
we obtain the following.

\begin{prop}
\label{pro:PN-corresp-G2} There is a sequence of pebbling moves which
places a pebble on $r$ in $G_{2}$ starting from $p_{2,f}$ if and
only if there is a setting of the existentially quantified variables
in $\phi$, which, together with $f$, satisfies $\phi$. 
\end{prop}
One useful property of the underlying graph together with an underlying
distribution is that it is not possible to accumulate more than five
pebbles on any vertex. This property will be instrumental in arguing
that the gadgets we attach to the underlying graph behave correctly.

\begin{prop}
\label{pro:underlying-graph-limit-5}It is not possible to place more
than five pebbles on any vertex in $G_{2}$ starting from any $p_{2,f}$.
\end{prop}
\begin{proof}
Suppose for a contradiction it is possible to place at least six pebbles
on a vertex $u$ in $G_{2}$. First, suppose $u$ is a vertex introduced
in our construction of $G_{2}$ from $G_{1}$, so that $u$ is an
internal vertex $w_{i}$, $1\leq i\leq\alpha$, in a one use path
$P=w_{0}w_{1}\ldots w_{\alpha}w_{\alpha+1}$. Let $D$ be a signature
of a minimum sequence of pebbling moves which places at least six
pebbles on $v_{i}$. Because $p_{2,f}(w_{i})=1$, for $\balance(D,p_{2,f},w_{i})\geq6$
we must have that the indegree of $w_{i}$ is at least five. It follows
by the pigeonhole principle that either the multiplicity of $w_{i-1}w_{i}$
or $w_{i+1}w_{i}$ is at least three. If the former is true, we can
apply \prettyref{cor:many-pebbles-begining-of-path} to obtain a sequence
of pebbling moves that places $2^{i}(3-1)+2\cdot1\geq6$ pebbles on
$w_{0}$. Similarly, if the latter is true, we apply \prettyref{cor:many-pebbles-begining-of-path}
to obtain a sequence of pebbling moves that places $2^{\alpha-i}(3-1)+2\cdot1\geq6$
pebbles on $w_{\alpha+1}$. Because $w_{0}$ and $w_{\alpha+1}$ are
vertices in $G_{1}$, it suffices to show that it is not possible
to place more than five pebbles on any vertex in $G_{1}$.

Suppose that $u$ is in $G_{1}$. Because it is possible to place
at least six pebbles on $u$ in $G_{2}$ starting from $p_{2,f}$,
by \prettyref{lem:one-use-edges}, there is a nonrepetitive sequence
of pebbling moves that places at least six pebbles on $u$ in $G_{1}$
starting from $p_{1,f}$. But this is clearly impossible, because
the maximum degree in $G_{1}$ is three and each vertex receives at
most two pebbles from $p_{1,f}$. 
\end{proof}
Now that we have established the important properties of the underlying
graph and the underlying distributions, we attach gadgets to the vertices
in the underlying graph. Just as the $\starg$ gadgets we attach in
our reduction from $\DPR$ to $\OPN$ force any potentially successful
distribution to take a certain form, our gadgets here force any potentially
unsuccessful distribution to take a form which effectively induces
one of the underlying distributions on the underlying graph.

\subsection{The Gadgets}

We introduce three classes of gadgets: the null gadget, the fork gadget,
and the eye gadget. In this section, we explore the relevant properties
of our gadgets as isolated graphs. 

All classes of gadgets share some common properties. The gadgets have
\emph{attachment vertices}; later, we will attach gadgets to the underlying
graph by identifying the attachment vertices of a gadget with vertices
in the underlying graph. A \emph{supply quota} $s$ assigns each attachment
vertex $v$ a number $s(v)$; each gadget has one or more supply quotas.
Under a particular distribution $q$, a gadget \emph{satisfies} $s$
if $q$ covers $s$.

The gadgets have \emph{overflow vertices,} which are adjacent to $r$;
we call the edges between the overflow vertices and $r$ the \emph{overflow
edges}. We say that a gadget has an \emph{overflow threshold} of $k$
if $r$ is reachable via an overflow edge under every distribution
of size $k$.

Let $q$ be a distribution of pebbles to a gadget. If the gadget is
able to satisfy any one of its supply quotas, or if $r$ is reachable
via an overflow edge, we say that the gadget is \emph{potent} under
$q$. We say that a gadget has a \emph{potency threshold} of $k$
if the gadget is potent under every distribution of $k$ pebbles.

Every gadget has one or more \emph{critical distributions}, each of
equal size. If $q$ is a critical distribution and $s$ is a supply
quota, we say that $q$ \emph{breaches} $s$ if there exists a vertex
$v$ such that it is possible to place more than $s(v)$ pebbles on
$v$ starting from $q$. 

Our critical distributions and supply quotas are in bijective correspondence;
that is, for each critical distribution there is a corresponding supply
quota and vice versa. Each critical distribution $q$ exhibits the
following \emph{critical distribution properties}: 

\begin{enumerate}
\item starting from $q$, $r$ is not reachable via an overflow edge
\item $q$ does not breach its corresponding supply quota
\end{enumerate}
As we present the gadgets, their supply quotas, and their critical
distributions, we will establish an overflow threshold, a potency
threshold, and the critical distribution properties.

To motivate the study of these parameters, we outline their use in
our proof of the correctness of our reduction. Given an $\RAETSAT$
instance $\phi$, we compute $H$ and $k$ such that $\phi$ is valid
if and only if $\pi(H,r)\leq k$. We construct $H$ by attaching various
gadgets to the underlying graph and we set $k$ to be the sum, over
all gadgets, of the size of the gadget's critical distributions. 

Suppose that $\phi$ is valid and consider a distribution of $k$
pebbles to $H$. If some gadget is assigned fewer pebbles than its
potency threshold, the pigeonhole principle implies that some gadget
receives more pebbles than its overflow threshold, and hence $r$
is reachable. Otherwise, all gadgets are potent. If $r$ is reachable
via some gadget's overflow edge, we are done. Otherwise, every gadget
is able to satisfy one of its supply quotas; this implies that our
initial distribution on $H$ covers some $p_{2,f}$. Because $\phi$
is satisfiable under $f$, we obtain from \prettyref{pro:PN-corresp-G2}
a sequence of pebbling moves in the underling graph which places a
pebble on $r$. 

The converse direction is somewhat trickier, but proceeds roughly
as follows. Suppose that $\pi(H,r)\leq k$ and consider a setting
$f$ of the universally quantified variables of $\phi$. We assign
pebbles to $H$ by selecting (according to $f$) a critical distribution
for each gadget. Because $\pi(H,r)\leq k$, we obtain a signature
$D$ of a minimum sequence of pebbling moves which places a pebble
on $r$. Next, we argue that our critical distribution properties
still apply even though the gadgets have been attached to the underlying
graph. Then we show how $D$ can be used to obtain a sequence of pebbling
moves in the underlying graph starting from $p_{2,f}$ which places
a pebble on $r$. A final application of \prettyref{pro:PN-corresp-G2}
implies that $\phi$ is satisfiable under $f$. 

Our gadgets are defined in terms of two parameters, $\beta$ and $c$.
We set $c=3$ (in fact, any constant $c$ so that $2^{c}$ exceeds
the constant obtained in \prettyref{pro:underlying-graph-limit-5}
will do). We postpone fixing the precise value of $\beta$; suffice
it to say we will choose $\beta=\Theta(\log n_{2})$. Our gadgets
use small paths of length $c$ to provide some separation between
the underlying graph and more sensitive areas of our gadgets. We use
larger paths of length $\beta$ so that the number of pebbles in a
gadget's critical distribution far exceeds its potency threshold.

\subsubsection{The Null Gadget}

\begin{figure}
\begin{center}\includegraphics[%
  scale=0.55]{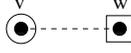}\end{center}

\caption{The null gadget. The dashed line represents a path of length $c$,
the circle around $v$ indicates that $v$ is an attachment vertex,
and the box around $w$ indicates that $w$ is an overflow vertex.
\label{fig:null-gadget}}
\end{figure}

The null gadget is a path of length $c$ and appears in \prettyref{fig:null-gadget}.
We use the null gadget to ensures that every vertex in the underlying
graph is not too far away from $r$, so that distributions which concentrate
pebbles on the underlying graph quickly imply that $r$ is reachable.
The null gadget has a single supply quota $s$, with $s(v)=0$; its
corresponding critical distribution $q$ assigns zero pebbles to each
vertex in the null gadget.

\begin{description}
\item [Overflow~threshold]Because $c$ is a fixed constant, the null gadget
is a fixed graph which does not depend upon $\phi$. By \prettyref{pro:PN-trivial-bound},
its pebbling number is a fixed constant, say $a$, not depending upon
$\phi$. Clearly, if there are $2a$ pebbles in the null gadget, then
it is possible to place two pebbles on $w$ and hence one pebble on
$r$. It follows that $2a=O(1)$ is an overflow threshold for the
null gadget. 
\item [Potency~threshold]Because $s$ is trivially satisfied, the null
gadget has a potency threshold of $0$. 
\item [Critical~distribution~properties]Because $q$ assigns zero pebbles
to the null gadget, it is clear that under $q$, the null gadget does
not breach $s$, nor is it possible to place a pebble on $r$ via
the null gadget's overflow edge.
\end{description}

\subsubsection{The Fork Gadget}

\begin{figure}
\begin{center}\includegraphics[%
  scale=0.55]{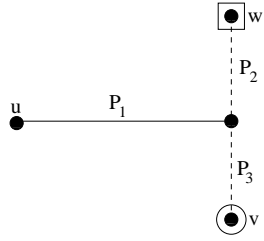}\end{center}

\caption{The fork gadget. The dashed lines represent paths $P_{2},P_{3}$
of length $c$, the solid line represents a path $P_{1}$ of length
$\beta$, the circle around $v$ indicates that $v$ is an attachment
vertex, and the box around $w$ indicates that $w$ is an overflow
vertex. \label{fig:fork-gadget}}
\end{figure}

The fork gadget consists of three paths $P_{1},P_{2},P_{3}$ which
share only a common endpoint, as shown in \prettyref{fig:fork-gadget}.
The fork gadget is responsible for injecting one pebble in the underlying
graph at the attachment location, much like the star gadgets in the
previous section. It has one supply quota $s$ with $s(v)=1$; the
corresponding critical distribution $q$ is given by $q(u)=2\cdot2^{\beta+c}-1$
and $q(x)=0$ for all $x\neq u$. 

\begin{description}
\item [Overflow~threshold]The fork gadget has an overflow threshold of
$2\cdot2^{\beta+c}+O(1)$. Indeed, if the fork gadget is unable to
place two pebbles on $w$ (and hence one on $r$), there can be at
most $O(1)$ pebbles on $P_{2}$ and $P_{3}$. Secondly, there can
be at most $2\cdot2^{\beta+c}-1$ pebbles in $P_{1}$ and $P_{2}$.
It follows that the fork gadget can contain at most $2\cdot2^{\beta+c}+O(1)$
pebbles if $r$ is not reachable via an overflow edge.
\item [Potency~threshold]The fork gadget has a potency threshold of $2^{\beta+c}+O(1)$.
Indeed, if the fork gadget is not potent, then it must have at most
$O(1)$ pebbles on $P_{2}$, or else it would be able to place a pebble
on $r$. Similarly, it must have at most $2^{\beta+c}-1$ pebbles
on $P_{1}$ and $P_{3}$, or else it would be able to place a pebble
on $v_{3}$ and therefore satisfy $s$.
\item [Critical~distribution~properties]Both the standard weight equation
and the greedy tree lemma show that under $q$, the fork gadget does
not breach $s$, nor is $r$ reachable via an overflow edge.
\end{description}

\subsubsection{The Eye Gadget}

\begin{figure}
\begin{center}\includegraphics[%
  scale=0.55]{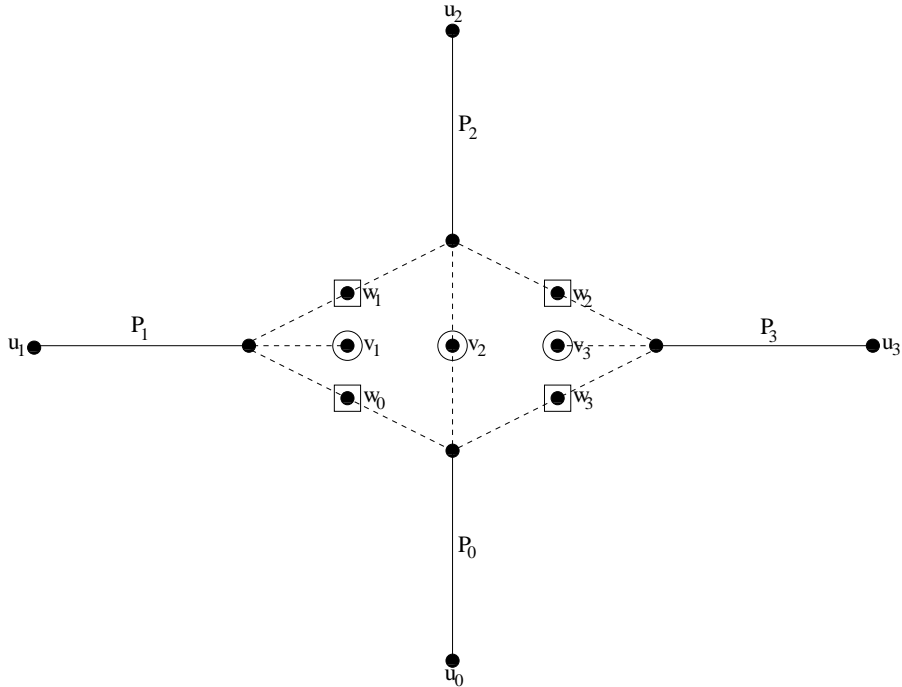}\end{center}

\caption{The eye gadget. The dashed lines represent paths of length $c$,
and the solid lines represent paths $P_{0},P_{1},P_{2},P_{3}$ of
length $\beta$. The circled vertices $v_{i}$ are attachment vertices;
the boxed vertices $w_{j}$ are overflow vertices. \label{fig:eye-gadget}}
\end{figure}

The eye gadget is the most complex of our three gadgets, and it is
at the heart of our reduction. Our reduction attaches one eye gadget
for each universally quantified variable in $\phi$. The eye gadget
is shown in \prettyref{fig:eye-gadget}. %

The eye gadget has two supply quota/critical distribution pairs. The
pair $(s^{+},q^{+})$ corresponds to a positive (true) setting of
the variable $x$ and the pair $(s^{-},q^{-})$ corresponds to a negative
(false) setting of $x$. We call $s^{+}$ the \emph{positive supply
quota} and we call $s^{-}$ the \emph{negative supply quota}. Similarly,
we call $q^{+}$ the \emph{positive critical distribution} and $q^{-}$
the \emph{negative critical distribution}. 

We define the supply quotas via $s^{+}(v_{1})=s^{+}(v_{3})=1$, $s^{+}(v_{2})=0$,
and $s^{-}(v_{1})=s^{-}(v_{3})=0$, $s^{-}(v_{2})=2$. Similarly,
the critical distributions are given by $q^{+}(u_{1})=q^{+}(u_{3})=2\cdot2^{\beta+c}-1$,
$q^{+}(u_{0})=q^{+}(u_{2})=2^{\beta+c}-1$, and $q^{-}(u_{1})=q^{-}(u_{3})=2^{\beta+c}-1$,
$q^{-}(u_{0})=q^{-}(u_{2})=2\cdot2^{\beta+c}-1$. 

Let $F$ be the subgraph of the eye gadget obtained by removing the
$u_{i}$ and all interior vertices of paths the $P_{i}$. Observe
that $F$ depends only on $c$ and therefore, like the null gadget,
$F$ is a fixed graph, not depending upon $\phi$. It follows that
$\pi(F)=O(1)$.

\begin{description}
\item [Overflow~threshold]The eye gadget has an overflow threshold of
$6\cdot2^{\beta+c}+O(1)$. Suppose the eye gadget contains $k$ pebbles
and it is not possible to place a pebble on $r$ via one of the overflow
edges. We show that $k\leq6\cdot2^{\beta+c}+O(1)$. Immediately, we
have that $F$ contains at most $2\pi(F)=O(1)$ pebbles, or else it
would be possible to place two pebbles on $w_{0}$ and hence one pebble
on $r$. To bound the number of pebbles in the $P_{i}$, we consider
two cases. First, suppose that each $P_{i}$ contains fewer than $2^{\beta+c}$
pebbles; in this case, we have that $k\leq4\cdot2^{\beta+c}+O(1)$.
Otherwise, suppose that $P_{j}$ has at least $2^{\beta+c}$ pebbles.
Clearly, $P_{j}$ has at most $2\cdot2^{\beta+c}-1$ pebbles, or else
we could use these pebbles to place a pebble on $r$ via the overflow
vertex $w_{j}$; similarly, the opposite path $P_{j+2}$ contains
at most $2\cdot2^{\beta+c}-1$ pebbles (subscript arithmetic is understood
modulo 4). Finally, the remaining paths $P_{j-1},P_{j+1}$ each contain
at most $2^{\beta+c}-1$ pebbles; indeed, if $P_{j-1}$ ($P_{j+1})$
contained $2^{\beta+c}$ pebbles, we could use them to place one pebble
on $w_{i-1}$ ($w_{i}$) and we could use $2^{\beta+c}$ pebbles from
$P_{j}$ to place a second pebble on $w_{i-1}$ ($w_{i}$). It follows
that the $P_{i}$ contain at most $6\cdot2^{\beta+c}-4$ pebbles,
and so $k\leq6\cdot2^{\beta+c}+O(1)$. 
\item [Potency~threshold]The eye gadget has a potency threshold of $5\cdot2^{\beta+c}+O(1)$.
Suppose the eye gadget contains $k$ pebbles, $r$ is not reachable
via an overflow edge, and it is not possible to satisfy $s^{+}$ or
$s^{-}$. We show that $k\leq5\cdot2^{\beta+c}+O(1)$. As before,
we have that $F$ contains at most $O(1)$ pebbles. To bound the number
of pebbles in the $P_{i}$, we consider the same two cases as before.
If each path has fewer than $2^{\beta+c}$ pebbles, we immediately
have $k\leq4\cdot2^{\beta+c}+O(1)$ and we're done. Otherwise, suppose
$P_{j}$ has at least $2^{\beta+c}$ pebbles. Once again, we have
that $P_{j}$ contains at most $2\cdot2^{\beta+c}-1$ pebbles, and
$P_{j-1},P_{j+1}$ each contain at most $2^{\beta+c}-1$. However,
now the opposite path $P_{j+2}$ has at most $2^{\beta+c}-1$ pebbles.
Indeed, if $P_{j},P_{j+2}$ both contain at least $2^{\beta+c}$ pebbles,
then we can either place one pebble each on $v_{1}$ and $v_{3}$,
satisfying $s^{+}$ (as is the case if $\left\{ j,j+2\right\} =\left\{ 1,3\right\} $),
or we can place two pebbles on $v_{2}$, satisfying $s^{-}$ (as is
the case if $\left\{ j,j+2\right\} =\left\{ 0,2\right\} $). It follows
that the paths contain at most $5\cdot2^{\beta+c}-4$ pebbles, implying
$k\leq5\cdot2^{\beta+c}+O(1)$.
\item [Critical~distribution~properties]It remains to verify the critical
distribution properties for $q^{+}$ and $q^{-}$. First, we show
that under $q\in\left\{ q^{+},q^{-}\right\} $, $r$ is not reachable
via an overflow vertex. Let $R$ be the set of overflow vertices in
the eye gadget, and let $D$ be the signature of a minimum sequence
of pebbling moves that places two pebbles on a vertex in $R$. By
\prettyref{lem:minimum-signatures-R}, we have that each vertex in
$R$ has outdegree zero in $D$. Observe that deleting $R$ from the
eye gadget results in a graph with three components; let $A_{1}$
be the component containing $P_{1}$, let $A_{2}$ be the component
containing $P_{0}$ and $P_{2}$, and let $A_{3}$ be the component
containing $P_{3}$. Let $T_{1}=A_{1}+\left\{ w_{0},w_{1}\right\} $,
let $T_{2}=A_{2}+R$, and let $T_{3}=A_{3}+\left\{ w_{2},w_{3}\right\} $,
as shown in \prettyref{fig:eye-trees}. %
\begin{figure}
\begin{center}\includegraphics[%
  scale=0.55]{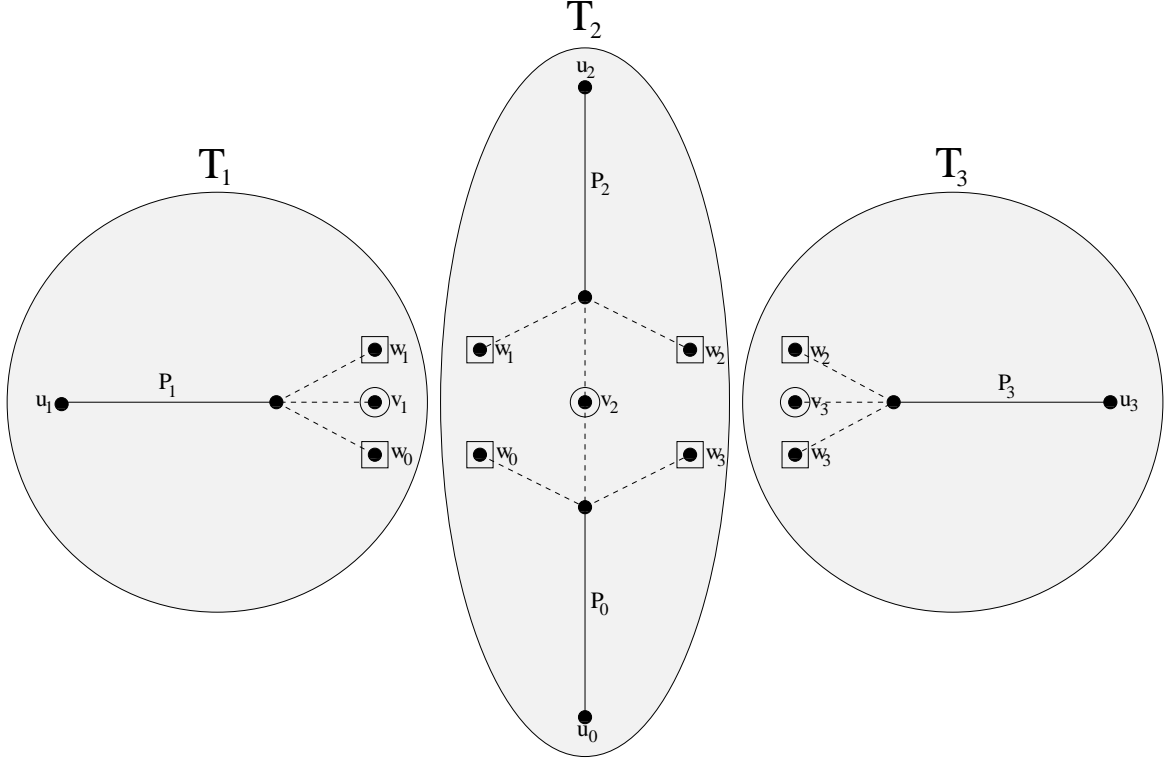}\end{center}

\caption{The overflow vertices split the eye gadget into three trees $T_{1},T_{2},T_{3}$\label{fig:eye-trees}}
\end{figure}
 Let $D_{l}$ be the digraph obtained by deleting from $D$ all pebbling
moves outside of $T_{l}$. Because $D$ is acyclic, it is immediate
that each $D_{l}$ is acyclic. Observe that for all $x\in V(T_{l})-R$,
we have $\balance(D_{l},q,x)=\balance(D,q,x)$. Furthermore, because
$d_{D}^{+}(w_{i})=0$ for all $w_{i}\in R$, we have that $d_{D_{l}}^{+}(w_{i})=0$.
Therefore by the acyclic orderability characterization, $D_{l}$ is
orderable. \\
Let $w_{i}$ be the overflow vertex on which $D$ places two pebbles.
Because $q(w_{i})=0$, we have that the indegree of $w_{i}$ in $D$
is at least two. Suppose two edges into $w_{i}$ are contained in
the same tree $T_{l}$. Then $D_{l}$ is the signature of a sequence
of pebbling moves in $T_{l}$ starting from $q$ that places at least
two pebbles on $w_{i}$. By the greedy tree lemma, the greedy pebbling
strategy in $T_{l}$ under $q$ places at least two pebbles on $w_{i}$.
However it is easily checked that regardless of $q\in\left\{ q^{+},q^{-}\right\} $,
$T_{l}\in\left\{ T_{1},T_{2},T_{3}\right\} $, and $w_{i}\in R$,
the greedy strategy in $T_{l}$ under $q$ places at most one pebble
on $w_{i}$. Alternatively, suppose that $D$ contains edges into
$w_{i}$ from two distinct trees. Because $w_{i}$ is in $T_{2}$
and one other tree, it must be that $D$ contains an edge into $w_{i}$
from $T_{2}$. Then $D_{2}$ is a signature of a sequence of pebbling
moves in $T_{2}$ starting from $q$ which places a pebble on $w_{i}$;
therefore the greedy strategy in $T_{2}$ starting from $q$ places
a pebble on $w_{i}$. Because the greedy strategy in $T_{2}$ starting
from $q^{+}$ is unable to place any pebbles on any overflow vertex,
it follows that $q=q^{-}$. Suppose that $D$ contains an edge into
$w_{i}$ from $T_{l}\in\left\{ T_{1},T_{3}\right\} $. Then $D_{l}$
is the signature of a sequence of pebbling moves in $T_{l}$ starting
from $q^{-}$ that places a pebble on $w_{i}$; therefore the greedy
strategy in $T_{l}$ starting from $q^{-}$ places a pebble on $w_{i}$.
But now a familiar contradiction is at hand: it is easily checked
that regardless of $T_{l}\in\left\{ T_{1},T_{3}\right\} $ and $w_{i}\in R$,
the greedy strategy in $T_{l}$ starting from $q^{-}$ is unable to
place a pebble on $w_{i}$.\\
Let $(s,q)\in\left\{ (s^{+},q^{+}),(s^{-},q^{-})\right\} $. It remains
to show that $q$ does not breach $s$. Suppose for a contradiction
that $D$ is the signature of a minimum sequence of pebbling moves
which witnesses that $q$ breaches $s$. We have that the outdegree
of each overflow vertex $w_{i}\in R$ is zero; indeed, if $d_{D}^{+}(w_{i})\geq1$,
then by \prettyref{lem:large-outdegree-many-pebbles} we would obtain
a sequence of pebbling moves placing two pebbles on $w_{i}$, a contradiction.
As before, let $D_{l}$ be the digraph obtained from $D$ by deleting
all edges outside of $T_{l}$; as before, we have that $D_{l}$ is
orderable in $T_{l}$. It follows that if $D$ places more than $s(v_{l})$
pebbles on $v_{l}$, then $D_{l}$ witnesses that it is possible to
place more than $s(v_{l})$ pebbles on $v_{l}$ in $T_{l}$ starting
from $q$. By the greedy tree lemma, the greedy strategy places more
than $s(v_{l})$ pebbles on $v_{l}$ in $T_{l}$ starting from $q$.
But now we have a contradiction: we easily check that regardless of
$(s,q)\in\left\{ (s^{+},q^{+}),(s^{-},q^{-})\right\} $ and $l\in\left\{ 1,2,3\right\} $,
the greedy strategy in $T_{l}$ starting from $q$ places exactly
$s(v_{l})$ pebbles on $v_{l}$. 
\end{description}

\subsubsection{Summary}

We summarize the various parameters of our gadgets in the following
table. 

\begin{center}\begin{tabular}{|c||c|c|c|}
\hline 
gadget&
potency threshold&
size of critical distributions&
overflow threshold\tabularnewline
\hline
\hline 
null&
$0$&
$0$&
$O(1)$\tabularnewline
\hline 
fork&
$2^{\beta+c}+O(1)$&
$2\cdot2^{\beta+c}-1$&
$2\cdot2^{\beta+c}+O(1)$\tabularnewline
\hline 
eye&
$5\cdot2^{\beta+c}+O(1)$&
$6\cdot2^{\beta+c}-4$&
$6\cdot2^{\beta+c}+O(1)$\tabularnewline
\hline
\end{tabular}\end{center}

From the table, we obtain the gap lemma.

\begin{lem}
[Gap Lemma]\label{lem:PN-gadget-gaps} There exists a nonnegative
constant $C$ (depending only on $c$) such that for each gadget,
the overflow threshold exceeds the size of the critical distributions
by at most $C$, and for the fork and eye gadgets, the size of the
critical distributions exceed the potency threshold by at least $2^{\beta+c}-C$. 
\end{lem}

\subsection{Construction of $H$. }

We set $\beta=\left\lceil \lg3Cn_{2}\right\rceil $, with $C$ as
in \prettyref{lem:PN-gadget-gaps}.

Armed with our gadgets and our underlying graph $G_{2}$, we are able
to describe the last step in our reduction from $\RAETSAT$ to $\RPN$.
For each pebble in $p_{2}$ on a vertex $z$ in the underlying graph,
we attach a fork gadget to $z$ by identifying the attachment vertex
$v$ in the fork gadget with $z$. For each triplet $z_{1},z_{2},z_{3}$
of vertices in $G_{2}$ corresponding to a universally quantified
variable $x$ in $\phi$, with $z_{1},z_{3}$ corresponding to positive
occurrences of $x$ in $\phi$ and $z_{2}$ corresponding to the negative
occurrence of $x$ in $\phi$, we attach an eye gadget by identifying
the attachment vertex $v_{i}$ in the eye gadget with $z_{i}$ in
the underlying graph. Finally, for any vertex $z\neq r$ in the underlying
graph to which we did not attach a fork or eye gadget, we attach a
null gadget by identifying $v$ in the null gadget with $z$ in the
underlying graph. Let $H$ be the resulting graph, and let $k$ be
the sum, over all gadgets in $H$, of the size of the gadget's critical
distributions. Our reduction from $\RAETSAT$ to $\RPN$ outputs $H$,
$k$, and $r$.

Note that we attach gadgets to the underlying graph by identifying
attachment vertices in gadgets with vertices in the underlying graph,
so that in $H$, each attachment vertex $v$ is a member of the underlying
graph and also a member of a gadget. Furthermore, by our construction,
every vertex other than $r$ in the underlying graph is identified
with an attachment vertex, so the vertices in the underlying graph
are exactly the attachment vertices together with $r$.

We pause to observe two important properties about $H$.

\begin{prop}
\label{pro:PN-num-gadgets-bound}In constructing $H$, we attach at
most two gadgets to every vertex in the underlying graph.
\end{prop}
\begin{proof}
Recall that $p_{2}$ assigns at most two pebbles to any vertex in
the underlying graph; furthermore, $p_{2}$ assigns at most one pebble
to any vertex associated with a universally quantified variable in
$\phi$. 
\end{proof}
\begin{prop}
\label{pro:diameter-H}The diameter of $H$ is at most $2\beta+O(1)$.
\end{prop}
\begin{proof}
It suffices to show that for each $z$ in $H$, the distance from
$z$ to $r$ is at most $\beta+O(1)$. If $z\neq r$, then $z$ is
contained in some gadget. In each gadget, every vertex is at most
$\beta+O(1)$ from an overflow vertex. 
\end{proof}

\subsection{$\RPNlong$ is $\ptc$-complete}

\begin{prop}
\label{pro:PN-pigenhole}Let $a_{1},\ldots,a_{n}$, $b_{1},\ldots,b_{n}$,
and $x$ be real numbers with $\sum_{i=1}^{n}a_{i}\geq\sum_{i=1}^{n}b_{i}$.
If $a_{n}<b_{n}-x$ then there exists $i$ such that $a_{i}>b_{i}+x/(n-1)$.
\end{prop}
\begin{proof}
By contradiction. Otherwise, \begin{eqnarray*}
\sum_{i=1}^{n}a_{i} & = & \sum_{i=1}^{n-1}a_{i}+a_{n}\\
 & < & \left(\sum_{i=1}^{n-1}b_{i}+\frac{x}{n-1}\right)+b_{n}-x\\
 & < & \sum_{i=1}^{n}b_{i}.\end{eqnarray*}

\end{proof}
We have accumulated the tools needed to show the correctness of our
reduction. 

\begin{thm}
\label{thm:phi-sat-iff-pi-H-leq-k}$\phi$ is valid if and only if
$\pi(H,r)\leq k$.
\end{thm}
\begin{proof}
($\implies$). Suppose that $\phi$ is valid and let $p$ be a pebble
distribution on $H$ of size $k$. We may assume $p(r)=0$. Let $l$
be the number of gadgets in $H$, label the gadgets as $Q_{1},\ldots,Q_{l}$,
let $a_{i}$ be the number of pebbles that $p$ assigns to $Q_{i}$,
and let $b_{i}$ be the size of $Q_{i}$'s critical distributions.
Because every vertex in $H$ besides $r$ belongs to at least one
gadget, we have $\sum_{i=1}^{n}a_{i}\geq k=\sum_{i=1}^{n}b_{i}$. 

We consider several cases. First, suppose there is some gadget $Q_{i}$
to which $p$ assigns fewer pebbles than $Q_{i}$'s potency threshold;
by the gap lemma, we have that $a_{i}<b_{i}-(2^{\beta+c}-C)$. By
\prettyref{pro:PN-pigenhole}, there is some $Q_{j}$ to which $p$
assigns at least $(2^{\beta+c}-C)/(l-1)$ pebbles more than $Q_{j}$'s
overflow threshold. By \prettyref{pro:PN-num-gadgets-bound}, $l-1\leq l\leq2n_{2}$.
It follows that $Q_{j}$ contains at least \begin{eqnarray*}
\frac{2^{\beta+c}-C}{2n_{2}} & \geq & \frac{2^{\beta}-C}{2n_{2}}\\
 & \geq & \frac{3Cn_{2}-C}{2n_{2}}\\
 & \geq & \frac{2Cn_{2}}{2n_{2}}\\
 & \geq & C\end{eqnarray*}
more pebbles than the size of its critical distributions. It follows
from \prettyref{lem:PN-gadget-gaps} that $Q_{j}$ contains at least
as many pebbles as its overflow threshold and therefore we can place
a pebble on $r$ via one of $Q_{j}$'s overflow edges. Otherwise,
$p$ assigns every gadget at least as many pebbles as its potency
threshold. If there is some gadget which is able to place a pebble
on $r$ via an overflow edge, then we are done. Otherwise, for every
gadget $Q$, there is a supply quota $s$ such that $Q$ under $p$
satisfies $s$. Using these supply quotas, we obtain a setting $f$
of the universally quantified variables in $\phi$ as follows. We
set $f(x)=\textrm{true}$ if the eye gadget associated with $x$ satisfies
its positive supply quota $s^{+}$; otherwise, the eye gadget associated
with $x$ must meet the negative supply quota $s^{-}$ and we set
$f(x)=\textrm{false}$. We claim that $p$ covers $p_{2,f}$. In each
gadget, execute the pebbling moves witnessing that the gadget satisfies
its supply quota. The fork gadgets alone produce a distribution that
is at least as good as $p_{2}$, and the eye gadgets supply the additional
pebbles proscribed by $p_{2,f}$. Because $\phi$ is valid, it follows
from \prettyref{pro:PN-corresp-G2} that $r$ is reachable.

($\Longleftarrow$). Suppose that $\pi(H,r)\leq k$ and let $f$ be
a setting of the universally quantified variables in $\phi$. We obtain
a setting of the existentially quantified variables in $\phi$ witnessing
that $\phi$ is satisfiable under $f$. Naturally, we study a pebble
distribution $p$ on $H$ of size $k$ corresponding to $f$; we construct
$p$ by choosing a critical distribution $q_{i}$ for each gadget
$Q_{i}$. If $Q_{i}$ is not an eye gadget, then $Q_{i}$ has only
one critical distribution and our selection of $q_{i}$ is forced.
If $Q_{i}$ is an eye gadget, we let $q_{i}$ be the positive critical
distribution $q^{+}$ if $f(x)=\textrm{true}$ and we let $q_{i}$
be the negative critical distribution $q^{-}$ otherwise. Note that
$p$ does not assign any pebbles to any vertex in the underlying graph.
Let $s_{i}$ be the supply quota associated with $q_{i}$. 

Let $H'$ be the graph obtained from $H$ by removing all the overflow
edges. Our first task is to establish the analog of \prettyref{pro:underlying-graph-limit-5}
for $H'$.
\begin{claim}
\label{cla:gadget-limit-5} In $H'$ starting from $p$, it is not
possible to place more than five pebbles on any vertex in the underlying
graph. 
\end{claim}
\begin{proof}
Suppose for a contradiction that $D$ is the signature of a minimum
sequence of pebbling moves that places at least six pebbles on some
vertex $w$ in the underlying graph. 
\begin{claim*}
$D$ does not contain an edge whose origin is inside the underlying
graph and whose destination is outside the underlying graph. 
\end{claim*}
\begin{proof}
Suppose for a contradiction that $uv$ is an edge in $D$ from a vertex
$u$ in the underlying graph to some vertex $v$ not in the underlying
graph. Because $H'$ does not contain any overflow edges, it must
be that $uv$ is an edge on a path of length $c$ in some gadget;
let this path be $P=x_{0}\ldots x_{c}$, with $u=x_{c}$ and $v=x_{c-1}$.
It follows that $D$ contains the edge $x_{1}x_{0}$, or else $D$
contains a cycle or a proper sink other than $w$, contradicting the
minimum signatures lemma. Because $p(x_{i})=0$ for each internal
vertex of $P$, we have by \prettyref{cor:many-pebbles-begining-of-path}
that it is possible to place $2^{c}=8\geq6$ pebbles on $u$ using
fewer pebbling moves, a contradiction. Therefore $D$ does not contain
an edge from the underlying graph to a vertex outside the underlying
graph.
\end{proof}
\begin{claim*}
For each $u$ in the underlying graph, the number of edges in $D$
into $u$ with origins outside the underlying graph is at most $p_{2,f}(u)$. 
\end{claim*}
\begin{proof}
If this were not the case, then there is some gadget $Q_{i}$ attached
to $u$ such that $D$ contains more than $s_{i}(u)$ edges from $Q_{i}$
into $u$. Construct $D'$ from $D$ by deleting all edges not contained
in $Q_{i}$. Clearly, $D'\subseteq D$ is acyclic; we show that $D'$
is orderable by verifying the balance condition. Consider a vertex
$v$ in $Q_{i}$. Recall that $H'$ does not contain overflow edges,
and therefore if $v$ is not an attachment vertex, then the neighborhood
of $v$ is contained in $Q_{i}$. It follows that if $v$ is not an
attachment vertex, we have $\balance(D',q_{i},v)=\balance(D,p,v)$.
Alternatively, if $v$ is an attachment vertex, we have that $d_{D'}^{+}(v)=0$,
or else $D'$ (and hence $D$) would contain an edge from a vertex
$v$ in the underlying graph to a vertex outside the underlying graph,
contradicting our previous claim. It follows that if $v$ is an attachment
vertex, we have $\balance(D',q_{i},v)\geq0$. By the acyclic orderability
characterization, we have that $D'$ is orderable under $q_{i}$.
Together with $d_{D'}^{-}(u)>s_{i}(u)$ and $d_{D'}^{+}(u)=0$ (recall
$u$ is an attachment vertex), we have that $\balance(D',q_{i},u)>s_{i}(u)$.
Therefore $D'$ witnesses that it is possible to place more than $s_{i}(u)$
pebbles on $u$ in $Q_{i}$ starting from $q_{i}$, contradicting
$Q_{i}$'s critical distribution properties. 
\end{proof}
We return to our proof of \prettyref{cla:gadget-limit-5}. Construct
$D'$ from $D$ by removing all edges from $D$ that are not in the
underlying graph. Clearly, $D'\subseteq D$ is acyclic. We show that
$D'$ is orderable under $p_{2,f}$ by checking the balance condition.
For each $u$ in the underlying graph, we have $\balance(D',p_{2,f},u)\geq\balance(D,p,u)$.
Indeed, at most $p_{2,f}(u)$ edges into $u$ are deleted from $D$
in our construction of $D'$; however, $p(u)=0$, so that $p_{2,f}$
offsets this decrease in balance. It follows that $D'$ is orderable
under $p_{2,f}$. Together with $\balance(D',p_{2,f},w)\geq\balance(D,p,w)\geq6$,
we have that it is possible to place at least six pebbles on $w$
starting from $p_{2,f}$ in the underlying graph, contradicting \prettyref{pro:underlying-graph-limit-5}.
This completes our proof of \prettyref{cla:gadget-limit-5}.
\end{proof}
We return to our proof of \prettyref{thm:phi-sat-iff-pi-H-leq-k}.
Let $D$ be the signature of a minimal sequence of pebbling moves
in $H$ starting from $p$ that places a pebble on $r$. 
\begin{claim}
[No Backflow into Gadgets Claim]\label{cla:no-backflow-into-gadgets}$D$
does not contain an edge from a vertex inside the underlying graph
to a vertex outside the underlying graph. 
\end{claim}
\begin{proof}
By the minimum signatures lemma, we have that $D$ contains at most
one pebbling move along an overflow edge and any such pebbling move
must be directed from an overflow vertex into $r$. Construct $D'$
from $D$ by removing this edge if it exists. Because $r$ has outdegree
zero in $D$, the acyclic orderability characterization implies that
$D'$ is orderable. Furthermore, because $D'$ does not contain any
pebbling move along overflow edges, $D'$ yields a sequence of pebbling
moves in $H'$. 

Because $D'$ is constructed from $D$ by removing at most one edge
into $r$, it suffices to show that $D'$ does not contain an edge
from a vertex inside the underlying graph to a vertex outside the
underlying graph. Suppose for a contradiction that $D'$ contains
an edge $uv$ from $u$ inside the underlying graph to $v$ outside
the underlying graph. It must be that $uv$ is a pebbling move along
a path $P$ of length $c$ in some gadget. Let $P=x_{0}\ldots x_{c}$
with $u=x_{c}$ and $v=x_{c-1}$. It follows that $D'$ contains the
edge $x_{1}x_{0}$. Indeed, if $D'$ does not have $x_{1}x_{0}$ as
an edge, neither does $D$ (after all, $x_{0}\neq r$), and so $D$
contains a cycle or a proper sink other than $r$, contradicting the
minimum signatures lemma. Therefore $D'$ contains the pebbling move
$x_{1}x_{0}$. 

Recalling that $p$ assigns each internal vertex of $P$ zero pebbles,
\prettyref{lem:large-outdegree-many-pebbles} implies that there is
an orderable $D''\subseteq D'$ which places at least $2^{c}=8$ pebbles
on $x_{c}=u$. But now $D''$ is a signature witnessing that it is
possible to place at least six pebbles on $u$ in $H'$ starting from
$p$, contradicting \prettyref{cla:gadget-limit-5}. 
\end{proof}
Let us resume our proof of \prettyref{thm:phi-sat-iff-pi-H-leq-k}.
Construct $D_{i}$ from $D$ by deleting from $D$ all edges not contained
in $Q_{i}$ or along $Q_{i}$'s overflow edges.
\begin{claim}
\label{cla:gadget-feasibility}$D_{i}$ is orderable under $q_{i}$,
and for each attachment vertex $v$, $\balance(D_{i},q_{i},v)=d_{D_{i}}^{-}(v)$.
\end{claim}
\begin{proof}
Because $D_{i}\subseteq D$, $D_{i}$ is acyclic and so it suffices
to verify the balance condition. Because $d_{D}^{+}(r)=0$, clearly
$d_{D_{i}}^{+}(r)=0$ and so the balance condition is satisfied at
$r$. Consider a vertex $v$ in $Q_{i}$. Unless $v$ is an attachment
vertex, all edges incident to $v$ in $D$ also appear in $D_{i}$,
and so $\balance(D_{i},q_{i},v)=\balance(D,p,v)$. Otherwise, if $v$
is an attachment vertex, then $d_{D_{i}}^{+}(v)=0$ or else $D$ would
contain an edge from a vertex in the underlying graph to a vertex
outside the underlying graph, contradicting \prettyref{cla:no-backflow-into-gadgets}.
Together with $q_{i}(v)=0$, it follows that $\balance(D_{i},q_{i},v)=d_{D_{i}}^{-}(v)$.
By the acyclic orderability characterization, $D_{i}$ is orderable
under $q_{i}$.
\end{proof}
\begin{claim}
\label{cla:underlying-graph-indegree-bound}For each $u$ in the underlying
graph, $D$ contains at most $p_{2,f}(u)$ edges from outside the
underlying graph into $u$. 
\end{claim}
\begin{proof}
Suppose that $u$ is a counterexample to the claim. If $u=r$, then
there is some gadget $Q_{i}$ such that $D$ contains an edge $wr$
into $r$ along one of $Q_{i}$'s overflow edges. But $D_{i}$ also
contains $wr$ and, by \prettyref{cla:gadget-feasibility}, $D_{i}$
is orderable under $q_{i}$. Clearly, $\balance(D_{i},q_{i},r)\geq1$
and therefore $r$ is reachable in $Q_{i}$ under $q_{i}$, contradicting
the critical distribution properties of $Q_{i}$. Otherwise, if $u\neq r$,
then there is some gadget $Q_{i}$ such that $D$ contains more than
$s_{i}(u)$ edges into $u$ from vertices in $Q_{i}$. But these edges
are also in $D_{i}$, so that $d_{D_{i}}^{-}(u)>s_{i}(u)$. By \prettyref{cla:gadget-feasibility},
$D_{i}$ is the signature of a sequence of pebbling moves in $Q_{i}$
under $q_{i}$ placing more than $s_{i}(u)$ pebbles on $u$, contradicting
$Q_{i}$'s critical distribution properties. 
\end{proof}
Let us complete our proof of \prettyref{thm:phi-sat-iff-pi-H-leq-k}.
Construct $E$ from $D$ by deleting from $D$ any edges outside the
underlying graph. We show that $E$ is orderable under $p_{2,f}$.
Clearly, $E\subseteq D$ is acyclic and therefore it suffices to check
the balance condition. Consider a vertex $u$ in the underlying graph,
and let $m$ be the number of edges into $u$ from outside the underlying
graph. In constructing $E$ from $D$, the balance of $u$ decreases
by $m$; by \prettyref{cla:underlying-graph-indegree-bound}, we have
$m\leq p_{2,f}(u)$. Because $p(u)=0$, changing distributions from
$p$ to $p_{2,f}$ increases the balance of $u$ by $p_{2,f}(u)$.
It follows that $\balance(E,p_{2,f},u)\geq\balance(D,p,u)$. Therefore
$E$ is orderable under $p_{2,f}$ and so $r$ is reachable in the
underlying graph under $p_{2,f}$. A final application of \prettyref{pro:PN-corresp-G2}
implies that $\phi$ is satisfiable under $f$. This completes our
proof of \prettyref{thm:phi-sat-iff-pi-H-leq-k}.
\end{proof}
We are now able to complete our proof that $\RPNlong$ is $\ptc$-complete. 

\begin{thm}
\label{thm:RPN-is-Pi_2-hard} $\RPNlong$ is $\ptc$-complete, even
when the diameter of $H$ is at most $O(\log n(H))$ and $k=\poly(n(H))$.
\end{thm}
\begin{proof}
We have already observed that $\RPN$ is in $\ptc$ and checked the
correctness of our reduction; it remains to check the diameter condition
on $H$ and that $H$ and $k$ are not too large relative to $\phi$
so that our reduction is computable in polynomial time. By \prettyref{pro:diameter-H},
the diameter of $H$ is at most $2\beta+O(1)=2\left\lceil \lg3Cn_{2}\right\rceil +O(1)$.
Because $n_{2}$ is the number of vertices in the underlying graph,
we have $n_{2}\leq n(H)$ and therefore the diameter of $H$ is at
most $2\left\lceil \lg3Cn(H)\right\rceil +O(1)=O(\log n(H))$. 

It remains to check the size condition on $H$ and $k$. Because $G_{1}$
has the same number of vertices as $G^{\NPR}(\phi)$, \prettyref{pro:GNPR-size}
implies that the size of $G_{1}$ is polynomial in the size of $\phi$.
Because the underlying graph $G_{2}$ is $\mathcal{S}(G_{1},\alpha)$
with $\alpha=\max\left\{ \lg2t,\,4\lg e(G_{1})\right\} $ and $t\leq2n(G_{1})$,
we have that the size of the underlying graph is polynomial in the
size of $G_{1}$. Observe that each gadget has size linear in $\beta=\Theta(\log n_{2})$.
Together with \prettyref{pro:PN-num-gadgets-bound}, we have that
the size of $H$ is polynomial in the size of $G_{2}$. It follows
that the size of $H$ is polynomial in the size of $\phi$. Finally,
every gadget's critical distribution size is at most $O(2^{\beta})=O(n_{2})$;
together with \prettyref{pro:PN-num-gadgets-bound}, we have that
$k$ is polynomial in $n_{2}$ and hence polynomial in $n(H)$.
\end{proof}

\subsection{$\PNlong$ is $\ptc$-complete.}

After having established \prettyref{thm:RPN-is-Pi_2-hard}, it is
relatively easy to show that $\PN$ is $\ptc$-complete. 

\begin{thm}
$\PNlong$ is $\ptc$-complete.
\end{thm}
\begin{proof}
We have already observed that $\PN$ is in $\ptc$. To how that $\PN$
is $\ptc$-hard, we reduce $\RPN$ to $\PN$. Let $G$ be a graph
with target vertex $r$, and let $k\geq0$ be an integer. We produce
$H$ and $k'$ so that $\pi(G,r)\leq k$ if and only $\pi(H)\leq k'$.
By \prettyref{thm:RPN-is-Pi_2-hard}, our reduction may assume that
the diameter $d$ of $G$ is at most $c'\lg n(G)$ for an absolute
constant $c'$ and $k=\poly(n(G))$. 

We construct $H$ and $k'$ as follows. Let $n=n(G)$ and set $\alpha=\left\lceil kn^{c'}\right\rceil $.
We let $H$ be the graph consisting of $\alpha$ copies of $G$ that
share $r$, so that $H-r$ is $\alpha$ disjoint copies of $G-r$.
We set $k'=\alpha k$. Observe that $k'$ and the size of $H$ are
polynomial in the size of $G$. It remains to show that $\pi(G,r)\le k$
if and only if $\pi(H)\leq k'$. 

($\implies$). Suppose $\pi(G,r)\leq k$. Consider a distribution
of $k'=\alpha k$ pebbles to $H$ and let $u$ be some target vertex
in $H$. Observe that $d(u,r)\leq d$ and therefore to place a pebble
on $u$, it suffices to show that we can place $2^{d}$ pebbles on
$r$. Our strategy is as follows. If there is some copy of $G$ with
at least $k$ pebbles, then we arbitrarily select a set $S$ of $k$
pebbles from this copy of $G$; because $\pi(G,r)\leq k$, we can
use these pebbles to place a pebble on $r$. We repeat this strategy
until we are unable to find a copy of $G$ with at least $k$ pebbles.
Let $s$ be the number of pebbles we are able to place on $r$ via
this strategy. Observe that after executing this strategy $s$ times,
at least $k\alpha-ks$ unused pebbles remain in $H$, and furthermore,
if more than $\alpha(k-1)$ unused pebbles remain in $H$, then some
copy of $G$ contains at least $k$ unused pebbles. It follows that
\begin{eqnarray*}
k\alpha-ks & \leq & \alpha(k-1)\end{eqnarray*}
 and therefore $s\geq\alpha/k\geq n^{c'}\geq2^{c'\lg n}\geq2^{d}$. 

($\Longleftarrow$). Suppose $\pi(H)\leq k'$, and let $p$ be a distribution
of $k$ pebbles to $G$. Naturally, we define a distribution $q$
of $k'=\alpha k$ pebbles to $H$ by distributing $k$ pebbles in
each copy of $G$ according to $p$. Let $D$ be the signature of
a minimum sequence of pebbling moves that places a pebble on $r$.
By the minimum signatures lemma, all edges of $D$ are contained in
a single copy of $G$. It follows that $r$ is reachable in $G$ under
$p$. 
\end{proof}

\section{Conclusions}

As we have seen, many graph pebbling problems on unrestricted graphs
are computationally difficult. We have seen that $\PRlong$ and $\OPNlong$
are both NP-complete. The authors believe it more likely than not
that $\PRlong$ remains NP-complete even when the graphs are restricted
to be planar. However, we have more hope that $\PRlong$ may fall
to P when the graphs are restricted to be outerplanar. It may be interesting
to investigate the computational complexity of these problems when
the inputs are restricted to be planar or outerplanar.

We have also seen that $\PNlong$ is $\ptc$-complete, and therefore
both NP-hard and coNP-hard. It follows that unless the polynomial
hierarchy collapses to the first level, $\PNlong$ is in neither NP
nor coNP. Consequently, given $G$ and $k$, it is unlikely that we
can compute in polynomial time a collection $\mathcal{P}$ of candidate
distributions of size $k$ such that if $\pi(G)>k$, then some vertex
in $G$ is not reachable from some $p\in\mathcal{P}$ (or else $\PN$
would be in NP).

We have shown that $\PClong$ and $\PRlong$ are both NP-complete;
however, the computational complexity of these problems diverges when
we introduce a universal quantifier over pebble distributions. When
we add such a quantifier to $\PClong$, we obtain the problem of determining
if $\gamma(G)\leq k$, which is possible in polynomial time \cite{cover-thm-1,cover-thm-2}.
The computational difficulties in $\PClong$ are smoothed out by the
consideration of all pebble distributions of size $k$: there is a
nice structure to the maximum pebble distributions from which a graph
cannot be covered with pebbles. On the other hand, by adding a universal
quantifier over all pebble distributions of size $k$ to $\PRlong$,
we obtain $\RPN$, which asks us to decide if $\pi(G,r)\leq k$. Instead
of observing a decrease in the computational complexity, we have stumbled
upon a $\ptc$-complete problem.

We recall that the graph pebbling community has shown a fair deal
of interest in developing necessary conditions and sufficient conditions
for equality in $\pi(G)=n(G)$. Of course, the ultimate goal is to
develop a characterization for when equality holds. We should remark
that our hardness result for $\PNlong$ does not suggest that any
such characterization need be complex from a computational point of
view. Indeed, our $\PNlong$ hardness result produces $G$ and $k$
with $k>n(G)$. It may be interesting to explore the complexity of
deciding whether $\pi(G)=n(G)$.

\paragraph*{Acknowledgements.}

We thank David Bunde, Jeff Erickson, and Sariel Har-Peled for helpful
suggestions throughout the revision process.


\begin{thebibliography}{CCFHPST}
\bibitem[BCCMW]{shameless-plug}D.P. Bunde, E.W. Chambers, D. Cranston, K. Milans, D.B. West, \emph{Pebbling
and Optimal Pebbling in Graphs}, preprint, 2004.
\bibitem[C]{history}F.R.K. Chung, \emph{Pebbling in hypercubes}, SIAM J. Disc. Math, Volume
2 (1989), 467-472.
\bibitem[CHH]{diameter-two-graphs}T.A. Clarke, R.A. Hochberg, and G.H. Hurlbert, \emph{Pebbling in diameter
two graphs and products of paths}, J. Graph Th. 25 (1997), 119-128.
\bibitem[CCFHPST]{cover-pebbling}B. Crull, T. Cundiff, P. Feltman, G.H. Hurlbert, L. Pudwell, Z. Szaniszlo,
and Z. Tuza, \emph{The Cover Pebbling Number of Graphs,} submitted
to Discrete Mathematics, 2004.
\bibitem[CHKT]{large-connectivity}A. Czygrinow, G.H. Hurlbert, H. Kierstead, and W.T. Trotter, \emph{A
note on graph pebbling}, Graphs and Combin., 18 (2002) 2, 219-225.
\bibitem[H99]{survey-paper}G.H. Hurlbert, \emph{A survey of graph pebbling}, Congressus Numerantium
139 (1999), 41-64.
\bibitem[M]{Moews}D. Moews, \emph{Pebbling graphs}, J. Combin. Theory (B) 55 (1992),
244-252.
\bibitem[M-web]{Moews-Online}D. Moews, \emph{Pebbling graphs}, http://xraysgi.ims.uconn.edu/dmoews/pebbling-graphs.ps
\bibitem[PSV]{PSV}L. Pachter, H.S. Snevily, B. Voxman, \emph{On pebbling graphs}, Congressus
Numerantium 107 (1995) 65-80.
\bibitem[P]{complexity}C.H. Papadimitriou, \emph{Computational Complexity}, Addison--Wesley,
1994.
\bibitem[S]{cover-thm-2}J. Sj\"ostrand, \emph{The Cover Pebbling Theorem}, preprint, http://www.arxiv.org/abs/math.CO/0410129,
2004.
\bibitem[VW]{cover-thm-1}A. Vuong, M. I. Wyckoff, \emph{Conditions for Weighted Cover Pebbling
of Graphs}, preprint, http://www.arxiv.org/abs/math.CO/0410410, 2004.
\bibitem[W]{watson}N.G. Watson, \emph{The Complexity of Pebbling and Cover Pebbling},
preprint, http://arxiv.org/abs/math.CO/0503511, 2005.
\end{thebibliography}
\end{document}